\documentclass{article}
\usepackage[russian]{babel}
\usepackage{amsmath,amsthm,amssymb}
\usepackage{bbm}
\usepackage{graphics,graphicx}

\newtheorem{theorem}{Теорема}
\newtheorem{lemma}{Лемма}
\newtheorem{prop}{Утверждение}
\newtheorem{corollary}{Следствие}
\theoremstyle{definition}
\newtheorem{definition}{Определение}
\newtheorem{example}{Пример}
\newtheorem{rema}{Замечание}

\newcommand{\V}{\mathcal V}
\newcommand{\Z}{\mathbb Z}
\newcommand{\N}{\mathbb N}
\newcommand{\R}{\mathbb R}
\newcommand{\psz}{($\Psi0$)}
\newcommand{\pso}{($\Psi1$)}
\newcommand{\pst}{($\Psi2$)}
\newcommand{\psth}{($\Psi3$)}

\begin{document}

\title{Слабые четности и функториальные отображения} 
\author{И.\,М.~Никонов}
\date{}
\maketitle

  \begin{abstract}
В настоящей работе рассматриваются объекты, имеющие два эквивалентных описания: как функториальные отображения и как слабые четности. Функториальные отображения позволяют преобразовывать узлы и продолжать посредством этого инварианты узлов. Мы вводим понятие максимальной слабой четности и описываем ее для узлов на фиксированной замкнутой ориентированной поверхности. При помощи найденной слабой четности строится проекция из виртуальных узлов в классические.
  \end{abstract}

\section*{Введение}

\begin{figure}[h]
\centering
\includegraphics[height=2cm]{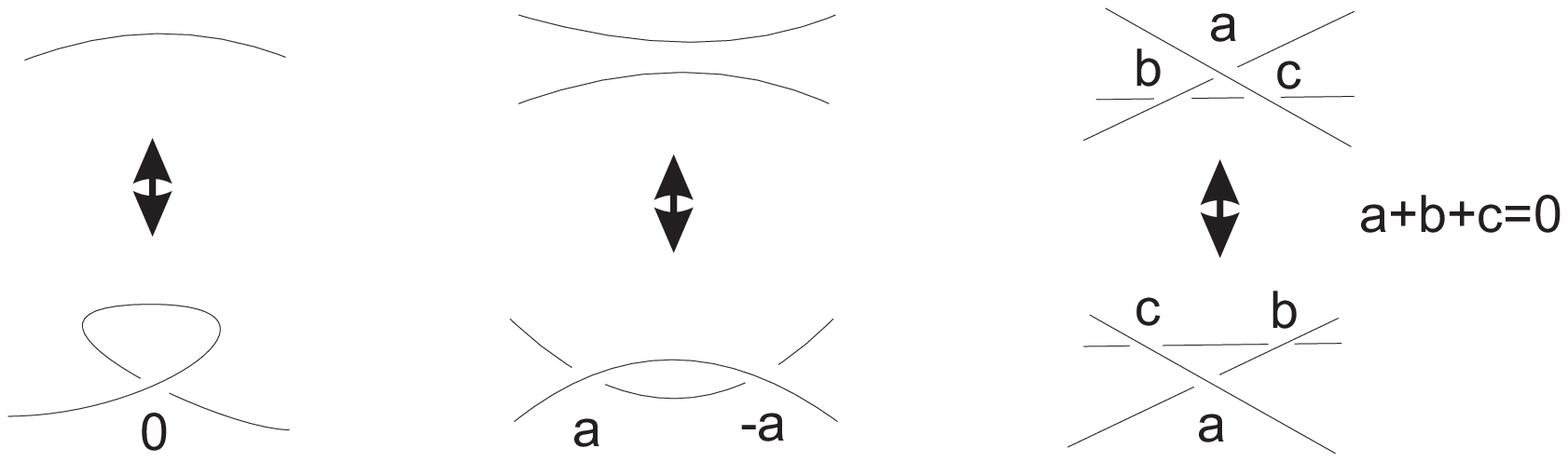}
\caption{Аксиомы четности.}\label{parity_axioms}
\end{figure}

Четность --- понятие, введенное В.О. Мантуровым в 2009
году~\cite{M1},--- определяется как способ оснащения перекрестков
диаграмм узла или зацепления числовыми метками $0$ и $1$, причем
приписывание меток должно быть согласовано с движениями
Рейдемейстера. Правило согласования можно сформулировать следующим
образом: четность перекрестков, не исчезающих и не появляющихся при
движении Рейдемейстера, не меняется; сумма (по модулю два) четностей
перекрестков, участвующих в движении Рейдемейстера, равна нулю (см.
рис.~\ref{parity_axioms}).

Введение дополнительной информации в диаграмму зацепления при помощи четностей позволяет усиливать инварианты узлов~\cite{M1,M2,M3,M12,M4,M5,M6,M7,IMN_mon}. При помощи новых инвариантов удалось показать нетривиальность свободных узлов, а также их классов кобордантности~\cite{M4}.

Другим приложением теории четностей является построение
функториальных отображений. Разделение перекрестков диаграммы на
четные и нечетные позволяет трансформировать диаграмму, по-разному
преобразуя ее в окрестности перекрестков различных типов, а именно,
заменять нечетные перекрестки на виртуальные, а четные оставлять как
есть. Из свойств четности   (на самом деле, из более слабых условий
{\em слабой четности}) следует, что результатом применения этого
преобразования к различным диаграммам одного узла будут диаграммы,
соответствующие некоторому другому, виртуальному узлу, то есть
функториальное отображение определено на уровне узлов, а не
диаграмм~\cite{M1}. Изучение функториальных отображений и слабых
четностей составляет содержание данной работы.

Статья организована следующим образом. В первом разделе даются основные определения. Мы начинаем с определения категории диаграмм узлов и приводим примеры некоторых таких категорий. Далее вводится определение функториального отображения и слабой четности. Затем обсуждаются простейшие свойства слабых четностей.

Второй раздел посвящен описанию слабых четностей для узлов в
утолщении фиксированной ориентированной замкнутой поверхности (а
также погруженных кривых в ориентированной замкнутой поверхности).
Мы доказываем, что гомотопическая слабая четность в данной теории
узлов является максимальной нетривиальной слабой четностью. Отсюда
следует, что для классических узлов все слабые четности тривиальны,
что обобщает аналогичный результат для четностей~\cite[Corollary
4.2]{IMN}. Завершает раздел описание проекции множества виртуальных
узлов на классические.


Автор выражает благодарность Д.\,П.~Ильютко и В.\,О.~Мантурову за
полезные обсуждения и комментарии.

Работа выполнена при частичной финансовой поддержке гранта Правительства РФ по постановлению N 220,
 договор 11.G34.31.0053,
РФФИ (грант \No~12-01-31432-а\_мол), Программы поддержки ведущих
научных школ РФ (грант \No~НШ-1410.2012.1), Программы Научные и
научно-педагогические кадры инновационной России
(госконтракт~14.740.11.0794).

\section{Функториальные отображения и слабые четности}

\subsection{Узлы и диаграммы}

Общепринятой формой представления узлов и зацеплений являются их диаграммы. Узлу соответствует бесконечно много различных диаграмм, но любые из них можно связать цепочкой преобразований --- изотопий диаграмм и движений Рейдемейстера.

\begin{definition}
Пусть $\mathcal K$ --- некоторый узел. {\em Категорией диаграмм} узла $\mathcal K$ называется категория $\mathfrak K$, объектами которой являются всевозможные диаграммы узла,  а морфизмами служат формальные композиции {\em элементарных морфизмов} --- изотопий диаграмм и движений Рейдемейстера.
\end{definition}

\begin{rema}
Наряду с категориями диаграмм конкретного узла можно рассматривать категории, состоящие из диаграмм узлов некоторого класса (классических узлов, виртуальных узлов, плоских узлов, свободных узлов, узлов и кривых на фиксированной поверхности и т.п.). Рассмотрение таких больших категорий удобно, когда изотопический класс узла заранее не определен, что имеет место, например, для образов функториальных отображений. Мы будем называть эти категории {\em теориями узлов}.
\end{rema}

Рассмотрим конкретные примеры типов узлов и их диаграмм.

\begin{definition}
{\em $4$-графом} называется любое несвязное объединение
четырехвалентных графов и {\em тривиальных компонент} ---
окружностей, каждая из которых рассматриваются как граф без вершин и
с одним (замкнутым) ребром. {\em Виртуальной диаграммой} называется
вложение в плоскость $4$-графа, вершины которого разделены на два
типа: {\em классические} и  {\em виртуальные}. В классических
вершинах выделена пара противоположных (полу)ребер, называемая {\em
проходом}; другие два (полу)ребра вершины образуют {\em переход}.
Вершины диаграммы называют также {\em перекрестками}. Диаграмма, не
содержащая виртуальных перекрестков, называется {\em классической}.
\end{definition}

\begin{rema}
При изображении виртуальной диаграммы виртуальные вершины обычно обводятся кружком. Проход в классической вершине обозначается разрывной линией, а переход --- сплошной (см. рис.~\ref{pic:virt_class_crossing})
\begin{figure}[h]
\centering
  \includegraphics[height=1cm]{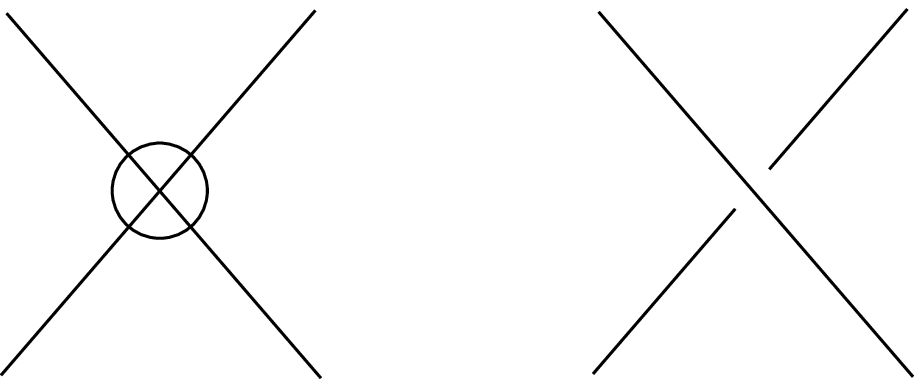}\\
  \caption{Виртуальный (слева) и классический (справа) перекрестки диаграммы }\label{pic:virt_class_crossing}
\end{figure}
\end{rema}

\begin {figure}[h]
\centering
\includegraphics[width=5cm]{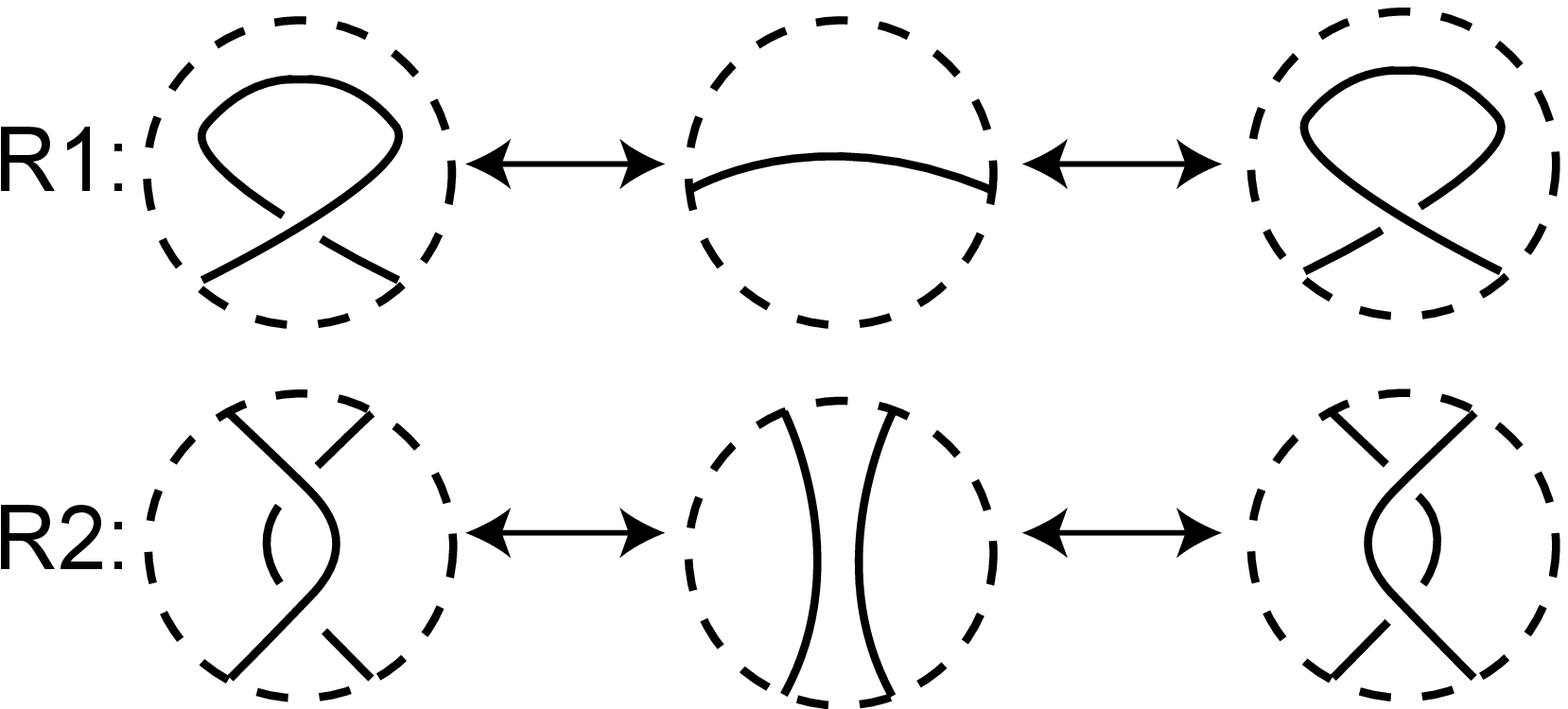}\qquad
\includegraphics[width=5cm]{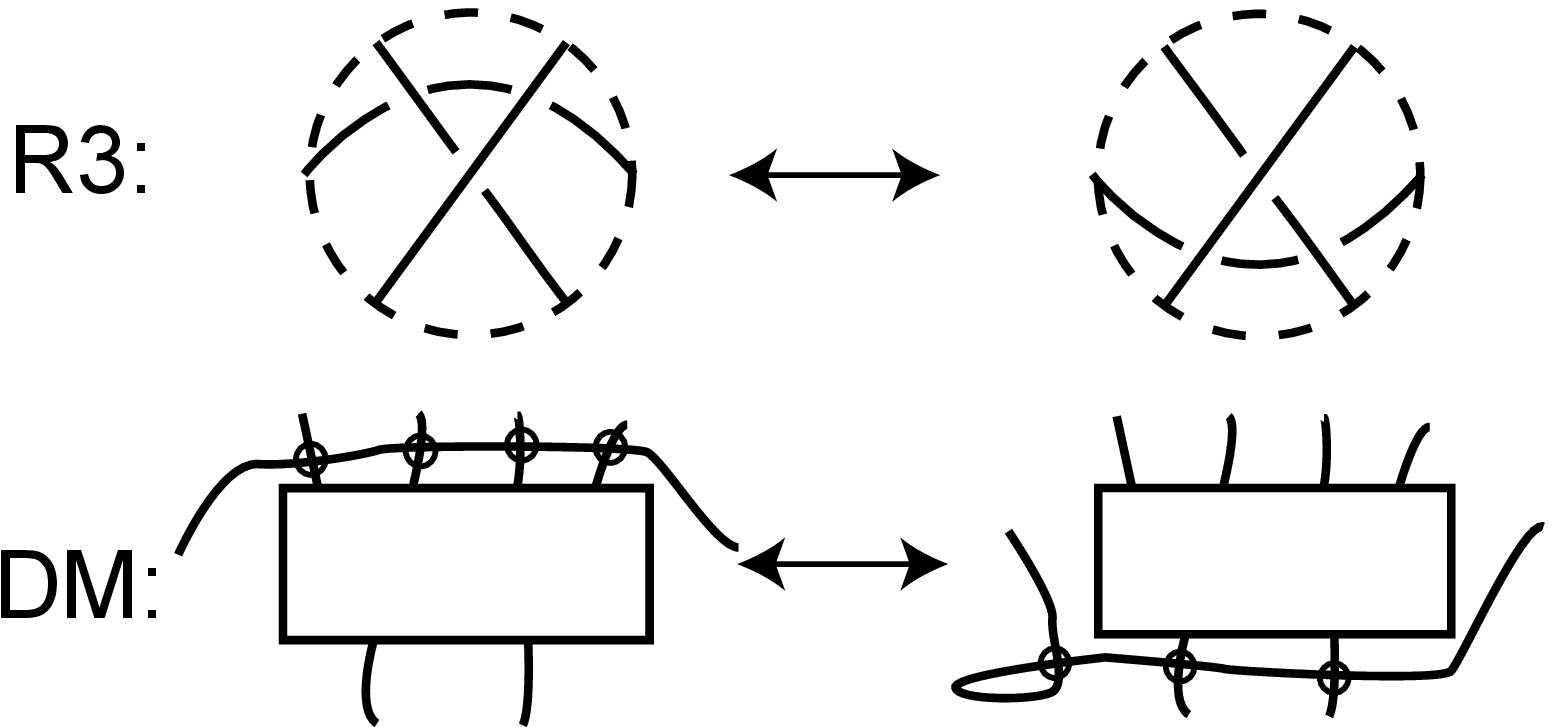}
\caption{Движения виртуальных диаграмм}\label{pic:virt_moves}
\end {figure}

\begin{definition}
{\em Движения виртуальных диаграмм} включают в себя обычные движения Рейдемейстера ($R1,R2,R3$), а также {\em движение объезда ($DM$)}, при котором произвольная дуга диаграммы, содержащая только виртуальные перекрестки, заменяется на новую дугу с теми же концами, также содержащую только виртуальные перекрестки (см. рис.~\ref{pic:virt_moves}). Класс эквивалентности виртуальных диаграмм по модулю движений называется {\em виртуальным зацеплением}.
\end{definition}

\begin{definition}
{\em Уникурсальной компонентой} виртуальной диаграммы называется
наименьшее по включению непустое множество ребер диаграммы,
замкнутое относительно переходов от ребра к противоположному ему (в
одном из его концов) ребру. Диаграмма с одной уникурсальной
компонентой называется {\em диаграммой виртуального узла}.
\end{definition}

\begin{rema}
Поскольку количество уникурсальных компонент не меняется при движениях диаграммы, корректно определено понятие {\em виртуального узла} как класса эквивалентности диаграмм с одной уникурсальной компонентой.
\end{rema}

\begin{definition}
Узел (зацепление), имеющий классическую диаграмму, называется {\em классическим}.
\end{definition}

\begin{rema}
Как следует из результата Куперберга~\cite{Kuperberg}, две классические диаграммы, соответствующие одному и тому же (классическому) зацеплению, можно связать последовательностью движений, в которой все промежуточные диаграммы тоже будут классическими. В силу этого, для классических узлов по умолчанию мы будем рассматривать категорию, состоящую из их классических диаграмм.
\end{rema}

\begin{rema}
По каждой виртуальной диаграмме $D$ можно построить диаграмму
$\widetilde D$, содержащую только классические перекрестки и лежащую
на некоторой ориентированной поверхности $S(D)$~\cite{KK}. Для этого
в каждом классическом перекрестке диаграммы располагается крест, а в
каждом виртуальном --- пара непересекающихся лент (см.
рис.~\ref{pic:virt_surface}). Соединяя эти кресты и ленты, мы
получим поверхность с краем $S'(D)$, куда вдоль осей лент
отображается диаграмма $\widetilde D$, получаемая из $D$ разведением
виртуальных перекрестков. Заклеивая компоненты края поверхности
$S'(D)$ дисками, мы получим искомую замкнутую поверхность $S(D)$.

\begin{figure}[h]
\centering
  \includegraphics[height=3cm]{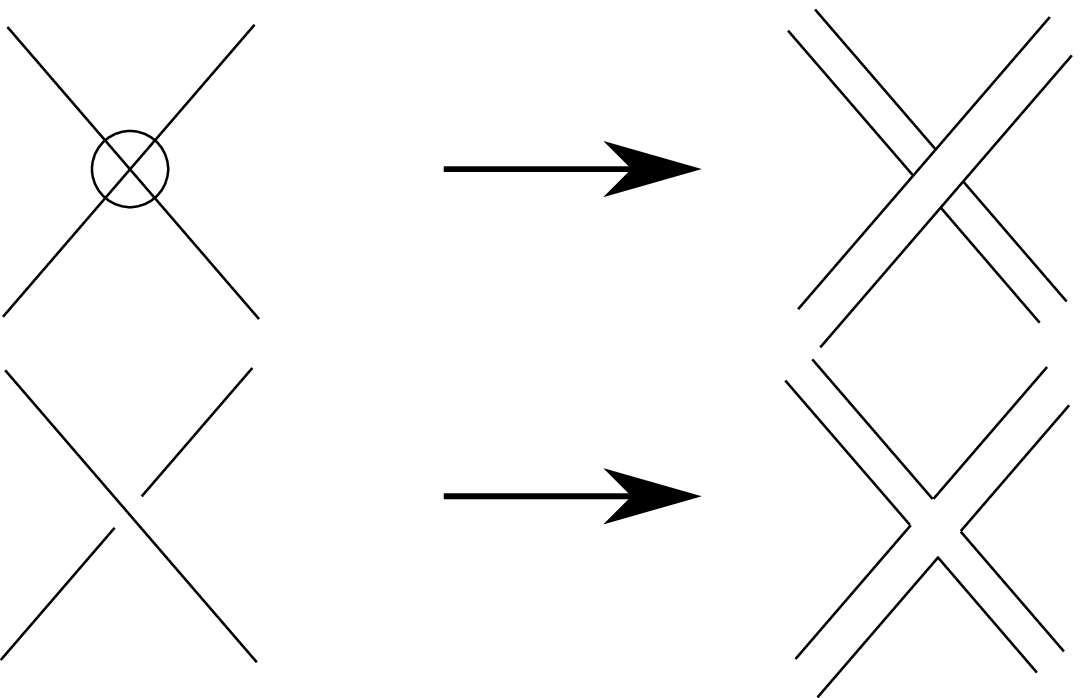}\\
  \caption{Построение поверхности $S'(D)$}\label{pic:virt_surface}
\end{figure}

Наоборот, по замкнутой ориентированной поверхности $S$  и диаграмме
$\widetilde D$ с классическими перекрестками в ней можно построить
виртуальную диаграмму следующим образом. Рассмотрим окрестность $S'$ диаграммы $\widetilde D$ в $S$.  Тогда $S'$ будет
ориентированной поверхностью с краем, так что имеется погружение
$S'$ в плоскость $\R^2$, сохраняющее ориентацию. При этом погружении
диаграмма $\widetilde D$ отобразится в некоторый $4$-граф $D$, в
котором образы перекрестков диаграммы $\widetilde D$ помечаются как
классические перекрестки (структура прохода-перехода при этом
индуцируется из диаграммы $\widetilde D$), а появляющиеся при
погружении точки самопересечения помечаются как виртуальные
перекрестки. Заметим, что диаграмма $D$ определена неоднозначно, но
оказывается, что различные погружения приводят к эквивалентным
виртуальным диаграммам.

Имеющееся соответствие между виртуальными диаграммами и диаграммами
на поверхностях позволяет дать новое определение виртуальных
зацеплений. Диаграмму на ориентированной поверхности $S$ можно
интерпретировать как проекцию вложения одной или нескольких
окружностей в утолщение поверхности $S\times[0,1]$. Отсюда выводится
(см.~\cite{KK}), что виртуальное зацепление можно определить как
вложение несвязного объединения окружностей в утолщения
ориентированных поверхностей, рассматриваемое  точностью до изотопии
и операции стабилизации (см. рис.~\ref{pic:stabilization}), меняющей
поверхность и ее утолщение.

\begin{figure}[h]
\centering
  \includegraphics[height=2cm]{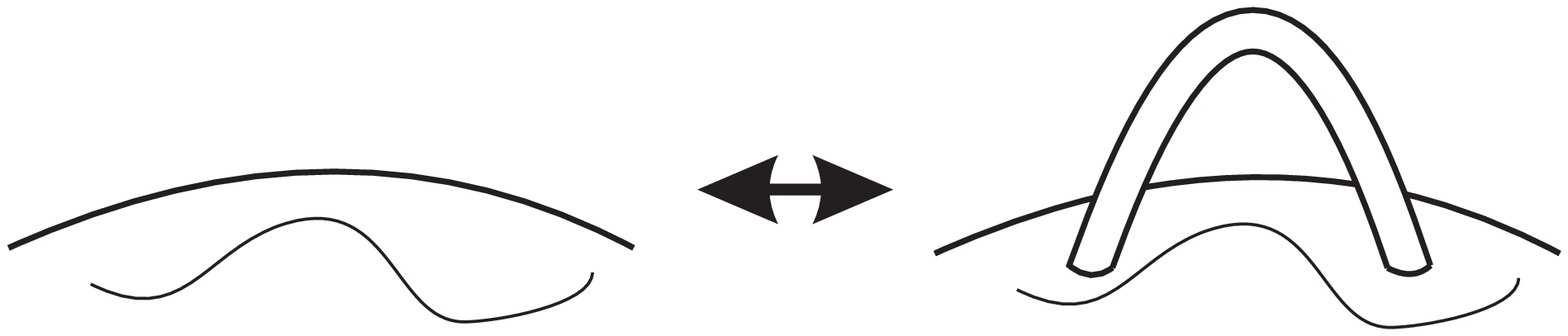}\\
  \caption{Операция стабилизации}\label{pic:stabilization}
\end{figure}

Последнее определение можно переформулировать следующим образом. Виртуальное зацепление --- это класс эквивалентности пар $(S,D)$, $S$ --- связная замкнутая ориентированная поверхность, $D$ --- диаграмма в $S$, у которой все перекрестки классические. Отношение эквивалентности порождено изотопиями, движениями Рейдемейстера и стабилизациями.
\end{rema}

Пусть $S$ --- двумерная связная замкнутая ориентированная поверхность.
\begin{definition}
{\em Зацеплением (соотв. узлом) на поверхности $S$} называется класс эквивалентности (классических) диаграмм (соотв. диаграмм с одной уникурсальной компонентой) на $S$  по модулю изотопий диаграмм и движений Рейдемейстера $R1,R2,R3$.
\end{definition}

\begin{rema}
1. Классические узлы и зацепления можно рассматривать как узлы на сфере $S^2$.

2. Имеется естественное отображение узлов на фиксированной поверхности $S$ во множество виртуальных узлов. На уровне диаграмм это отображение есть вложение категорий: диаграмме $D$ сопоставляется пара $(S,D)$.
\end{rema}

Забывание структуры проход-переход в классических перекрестках диаграмм приводит к новым теориям узлов.
\begin{definition}
{\em Плоским зацеплением} называется класс эквивалентности
виртуальных диаграмм по модулю движений и замен проходов на переходы
в любых перекрестках диаграммы. Классы эквивалентности виртуальных
диаграмм по модулю замен проходов на переходы в перекрестках
диаграммы называются {\em плоскими диаграммами}. Их можно задавать,
как вложения в плоскость $4$-графов, вершины которых разделены на
виртуальные и классические перекрестки, при этом структура
"проход-переход" в классических перекрестках не определяется. {\em
Плоским узлом} называется плоское зацепление с одной уникурсальной
компонентой.
\end{definition}

Для плоских узлов в качестве объектов категории диаграмм мы будем
рассматривать плоские диаграммы, а в качестве морфизмов --- движения
Рейдемейстера и объезда, индуцированные соответствующими движениями
на виртуальных диаграммах.

\begin{rema}
Забывание проходов и переходов в диаграммах узлов на фиксированной поверхности $S$ порождает множество вложений $4$-графов в поверхность, которые можно интерпретировать как погружения общего положения для набора окружностей. Классами эквивалентности диаграмм в этом случае будут классы свободной гомотопии набора кривых на поверхности.
\end{rema}

\subsection{Функториальные отображения}

\begin{definition} Пусть $\mathcal K$ --- некоторый узел.
{\em Функториальным отображением} на категории диаграмм  $\mathfrak
K$ узла $\mathcal K$ называется функтор $\Psi$ из категории
$\mathfrak K$ в категорию диаграмм виртуальных узлов (если узел
$\mathcal K$ есть виртуальный узел или узел на фиксированной
поверхности) либо категорию плоских диаграмм (если $\mathcal K$ ---
плоский узел или свободный гомотопический класс кривой на
фиксированной поверхности), такое что
\begin{enumerate}
\item для любой диаграммы $D\in ob(\mathfrak K)$ диаграмма $\Psi(D)$
получается из $D$ заменой некоторых классических перекрестков на виртуальные;
\item для любого движения Рейдемейстера $f\colon D_1\to D_2$
отображение $$\Psi(f)\colon \Psi(D_1)\to\Psi(D_2)$$ представляет
собой {\em то же} движение Рейдемейстера либо движение объезда.
\end{enumerate}
\end{definition}

Рассмотрим диаграмму $D\in ob(\mathfrak K)$. Обозначим через $\V(D)$
множество классических перекрестков диаграммы. Любой элементарный
морфизм $f\colon D\to D'$ (т.е. изотопия, движение Рейдемейстера или
движение объезда) определяет соответствие между перекрестками
диаграмм --- некоторое частично определенное отображение $f_*\colon
\V(D)\to \V(D')$. Областью определения (соотв., множеством значений)
отображения $f_*$, если $f$ --- изотопия, движение объезда или
третье движение Рейдемейстера, будет все множество $\V(D)$ (соотв.,
$\V(D')$) и, если $f$ --- первое или второе движения, будет
множество перекрестков диаграммы $D$ (соотв., диаграммы $D'$), не
участвующих в движении. Отображение $f_*$ представляет собой биекцию
области определения на множество значений.

Пусть $\Psi$ --- функториальное отображение. Тогда на множестве классических перекрестков $\V(D)$
каждой диаграммы $D\in ob(\mathfrak K)$ определено отображение $\psi_D\colon \V(D)\to
\Z_2$, такое что $\psi_D(v)=0$, если перекресток $v$ сохраняется в
диаграмме $\Psi(D)$, и $\psi_D(v)=1$, если перекресток $v$ в
диаграмме $\Psi(D)$ становится виртуальным. Мы будем называть перекресток $v$
{\em четным} относительно функториального отображения $\Psi$, если $\psi_D(v)=0$, и {\em нечетным}, если
$\psi_D(v)=1$.

\begin {figure}[h]
\centering
\includegraphics[width=10cm]{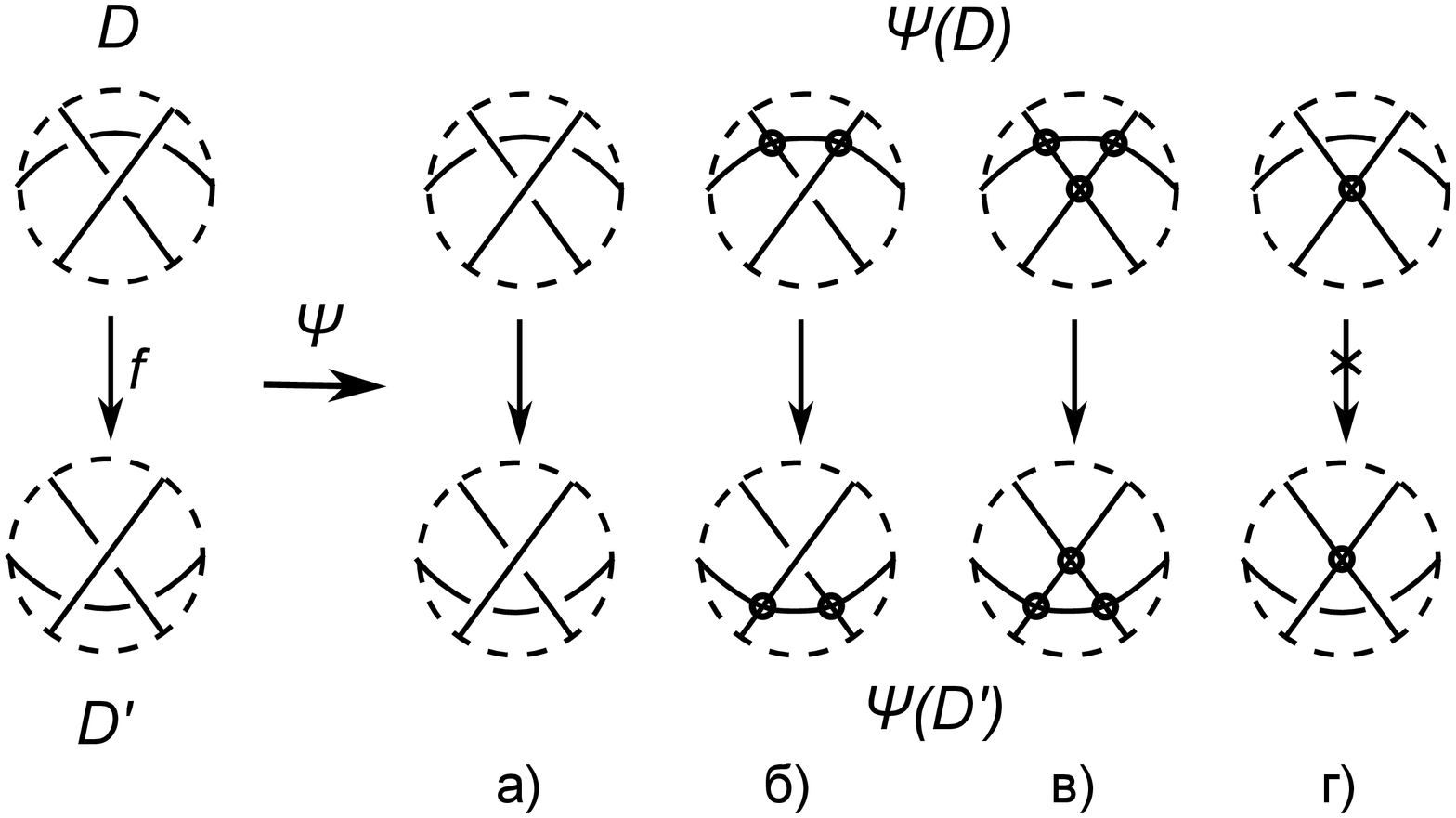}
\caption{Доказательство свойства ($\Psi3$): а) без замены
перекрестков остается третье движение Рейдемейстера; б,в) при замене
двух или трех перекрестков на виртуальные возникает движение
объезда; г) при замене одного перекрестка появляется запрещенное
движение.}\label{pic:psi3_prop}
\end {figure}

\begin{prop}
Отображения $\psi_D$ обладают следующими свойствами:
\begin{itemize}
\item[($\Psi0$)] для любого движения Рейдемейстера $f\colon D\to
D'$ и любого перекрестка $v\in\V(D)$, не участвующего в движении,
$\psi_{D}(v)=\psi_{D'}(f_*(v))$;
\item[($\Psi2$)] если $f\colon D\to D'$ --- уменьшающее второе
движение Рейдемейстера, в котором участвуют перекрестки $u,v\in\V(D)$,
то $\psi_D(u)=\psi_D(v)$;
\item[($\Psi3$)] если $f\colon D\to D'$ --- третье
движение Рейдемейстера, в котором участвуют перекрестки $u,v,w\in\V(D)$,
то $\psi_D(u)=\psi_{D'}(f_*(u))$, $\psi_D(v)=\psi_{D'}(f_*(v))$,
$\psi_D(w)=\psi_{D'}(f_*(w))$ и $\psi_D(u)+\psi_D(v)+\psi_D(w)\ne 1$
(как элемент $\Z$).
\end{itemize}
\end{prop}
\begin{proof}
Пусть перекресток $v$ не участвует в движении. Тогда соответствующий перекресток $\Psi(v)$ в диаграмме $\Psi(D)$ тоже не принимает участия в движении $\Psi(f)$. Следовательно, перекресток $\Psi(v)$ и его образ $\Psi(f_*(v))$ в диаграмме $\Psi(D')$ имеют одинаковый тип (классический либо виртуальный). Таким образом, $\psi_{D}(v)=\psi_{D'}(f_*(v))$, и свойство ($\Psi0$) выполнено.

Пусть $f$ --- третье движение Рейдемейстера  на перекрестках
$u,v,w$. Чтобы после замены некоторых перекрестков на виртуальные
получилось третье движение на тех же вершинах или движение объезда,
нужно сделать виртуальными $0,2$ или $3$ перекрестка из $u,v,w$;
замена ровно одного перекрестка приводит к запрещенному движению
(см. рис.~\ref{pic:psi3_prop}). Кроме того, четности у
соответствующих перекрестков диаграмм $\Psi(D)$ и $\Psi(D')$ должны
совпадать. Таким образом, верно свойство ($\Psi3$).

Свойство ($\Psi2$) доказывается аналогично.
\end{proof}

\begin{definition}
Cемейство отображений $\psi_D\colon \V(D)\to \Z_2$, $D\in
ob(\mathfrak K)$, для которого выполнены условия $(\Psi0), (\Psi2),
(\Psi3)$, называется {\em слабой четностью}.
\end{definition}

Заметим, что сопоставление функториальному отображению
слабой четности, описанное выше, устанавливает взаимно-однозначное
соответствие между этими классами. Ниже мы будем, в основном,
использовать понятие слабой четности.

\begin{example}[Нулевая слабая четность]
Определим слабую четность $\mathbb{O}$ при помощи равенства $\mathbb
O_D(v)\equiv 0$. Нулевой четности соответствует тождественное
функториальное отображение.
\end{example}

\begin{example}[Тривиальная слабая четность]
Тривиальная слабая четность $\mathbbm 1$ задается тождеством
$\mathbbm 1_D(v)\equiv 1$. Ей соответствует функториальное
отображение, которое заменяет все перекрестки на виртуальные. В
результате получаются чисто виртуальные диаграммы тривиального узла.
\end{example}

В дальнейшем мы будем изучать исключительно нетривиальные слабые
четности.

\begin{theorem}\label{th_psi1}
Пусть слабая четность $\psi$ нетривиальна. Тогда она обладает
свойством
\begin{itemize}
\item[($\Psi1$)] если $f\colon D\to D'$ --- уменьшающее первое
движение Рейдемейстера, в котором участвует вершина $u\in\V(D)$, то
$\psi_D(u)=0$.
\end{itemize}
\end{theorem}
\begin{proof}
Пусть $\psi$ --- нетривиальная слабая четность в категории диаграмм
$\mathfrak K$. Это значит, что найдется диаграмма $D$ и перекресток
в ней $v$, такой что $\psi_D(v)=0$. Применим к диаграмме $D$ второе
движение Рейдемейстера, как показано на рисунке~\ref{f_psi1}. По
свойству \pst имеем $\psi(w)=\psi(v)=0$. На дуге $vw$ добавим при
помощи первого движения Рейдемейстера перекресток $u$. Из условия
\psth, примененному к треугольнику $uvw$, следует, что $\psi(u)=0$.
Теперь уберем при помощи второго движения Рейдемейстера
вспомогательные перекрестки $w, w'$. Таким образом, перекресток $u$,
появляющийся при первом движении рядом с четным перекрестком $v$,
является четным. Обозначим получившуюся диаграмму через $D'$.

\begin{figure}
\centering
  \includegraphics[height=5cm]{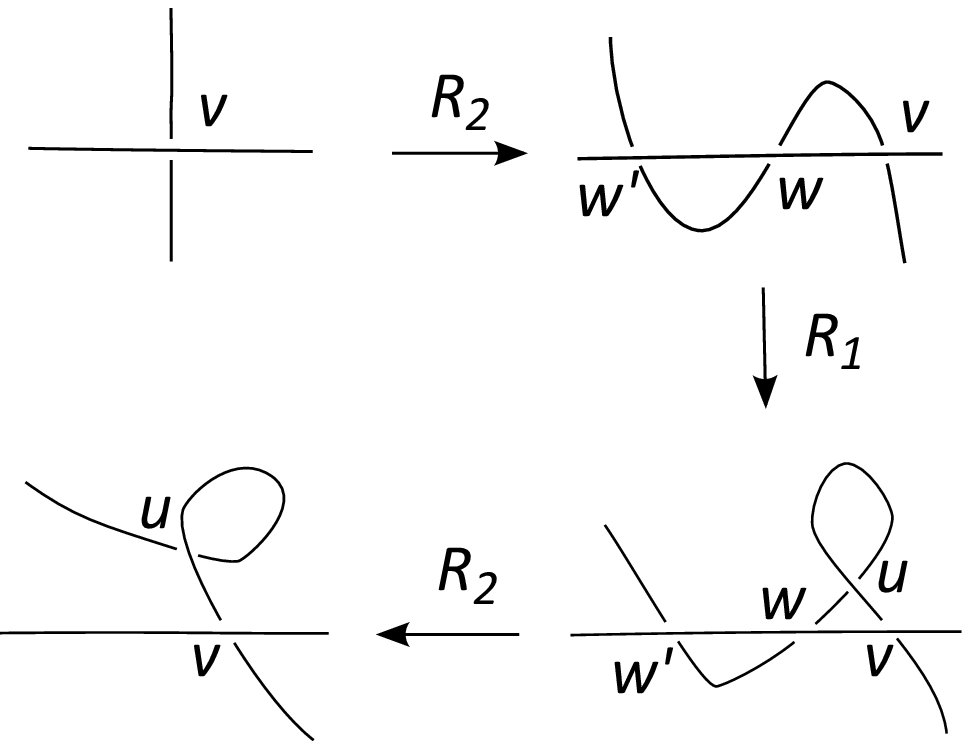}\\
  \caption{Доказательство свойства ($\Psi1$).}\label{f_psi1}
\end{figure}

Пусть $f\colon D_1\to D_2$ --- увеличивающее первое движение
Рейдемейстера в категории диаграмм $\mathfrak K$, и $t\in\V(D_2)$
--- новый перекресток. Так как диаграммы $D$ и $D_1$ относятся к
одному узлу, имеется последовательность движений Рейдемейстера,
переводящая $D$ в $D_1$. Тогда мы можем рассмотреть аналогичную
последовательность движений (добавляя при необходимости вторые и
третьи движения), которая переводит диаграмму $D'$ в $D_2$, причем
перекресток $u$ переходит в перекресток $t$. По свойствам \psz и
\psth получаем $\psi_{D_2}(t)=\psi_{D'}(u)=0$.
\end{proof}

Большое количество примеров нетривиальных слабых четностей можно
получить из обычных четностей.

Пусть $G$ --- некоторая группа. Предположим, что узел $\mathcal K$
ориентирован, так что на всех диаграммах категории диаграмм узла
$\mathfrak K$ задана индуцированная ориентация.

\begin{definition}
{\em Ориентированной четностью} $p$ на категории диаграмм $\mathfrak K$ называется семейство отображений $p_D\colon \V(D)\to G$, $D\in ob(\mathfrak K)$, обладающее следующими свойствами:
\begin{itemize}
\item для любого элементарного морфизма $f\colon D\to
D'$ и любого перекрестка $v\in\V(D)$ из области определения отображения $f_*$ выполнено равенство
$p_{D}(v)=p_{D'}(f_*(v))$;

\item если $f\colon D\to D'$ --- уменьшающее первое
движение Рейдемейстера и $v\in\V(D)$ --- исчезающий перекресток, то $p_D(v)=1$, где $1\in G$ --- единица группы;

\item
 если $f\colon D\to D'$ --- уменьшающее второе
или третье движение Рейдемейстера, то $\prod p_D(v_i)^{\epsilon(v_i)}=1$, где $v_i$  --- вершины, участвующие в движении, $\epsilon(v_i)$ --- индекс инцидентности вершины по отношению к исчезающей области (двуугольнику или треугольнику), см. рис.~\ref{pic:incidence_index} слева, а порядок вершин в произведении определяется ориентацией поверхности, в которой находится диаграмма.

\begin {figure}[h]
\centering
\includegraphics[height=2.5cm]{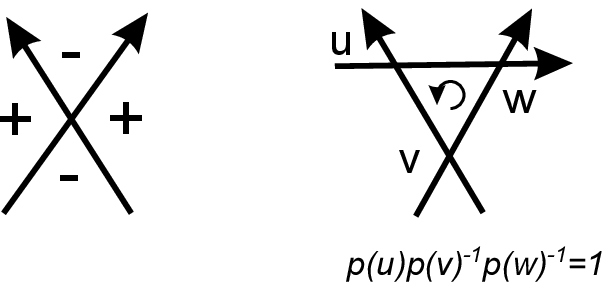}
\caption{Индекс инцидентности и соотношения ориентированной четности}\label{pic:incidence_index}
\end {figure}

\end{itemize}
\end{definition}

На рис.~\ref{pic:incidence_index} справа указан пример соотношения для ориентированной четности.

\begin{rema}
1. Для любой ориентированной четности $p$ и уменьшающего второго движения Рейдемейстера $f\colon D\to D'$ на перекрестках $u$ и $v$ выполняется соотношение
$$ p_D(u)p_D(v)=1,$$
поскольку индексы инцидентности перекрестков $u$ и $v$ совпадают.

2. Ориентированная четность не зависит от ориентации узла. Более точно, если семейство отображений $p$ является ориентированной четностью на категории диаграмм узла $\mathcal K$, то оно же является ориентированной четностью на категории диаграмм узла $-\mathcal K$, получающегося из $\mathcal K$ обращением ориентации (здесь мы отождествляем естественным образом категории диаграмм узлов $\mathcal K$ и $-\mathcal K$).

3. Четности с коэффициентами, определенные в работе~\cite{IMN}, являются ориентированными четностями в силу аксиом четности и соотношений вида $2p(v)=0$ (см.~\cite[Lemma 4.4]{IMN}).
\end{rema}

\begin{prop}
Пусть $p$ --- (ориентированная) четность на категории диаграмм
$\mathfrak K$. Тогда семейство отображений $\psi_D\colon
\V(D)\to\Z_2$, определяемое по правилу
$$
\psi_D(v)=\left\{
            \begin{array}{cl}
              1, & p_D(v)\ne 1;\\
              0, & p_D(v)= 1,
            \end{array}
          \right.
$$
является нетривиальной слабой четностью на категории диаграмм $\mathfrak K$. Эта
слабая четность называется {\em индуцированной слабой четностью}.
\end{prop}

\begin{proof}
Свойство ($\Psi0$) следует из первого условия в определении
ориентированной четности, а свойство ($\Psi1$) --- из второго
условия. Таким образом, $\psi$ не является тривиальной слабой
четностью. Из третьего условия определения ориентированной четности
вытекает, что четности перекрестков, участвующих во втором движении,
одновременно равны либо не равны единице группы коэффициентов, что
означает, что выполнено свойство ($\Psi2$). Аналогично, если два
перекрестка из трех, участвующих в третьем движении Рейдемейстера,
имеют тривиальную четность относительно $p$, то четность третьего
перекрестка также тривиальна. Это доказывает свойство ($\Psi3$).
\end{proof}

\subsection{Операции на слабых четностях}

На слабых четностях можно определить естественное отношение
частичного порядка.
\begin{definition}
Пусть $\psi,\psi'$ --- слабые четности в категории диаграмм $\mathfrak K$. Будем говорить, что
слабая четность $\psi$ {\em не превосходит} слабой четности $\psi'$ ($\psi\le\psi'$),
если для любой диаграммы $D\in ob(\mathfrak K)$ и любой вершины $v\in\V(D)$ выполнено
неравенство $\psi_D(v)\le\psi'_D(v)$.
\end{definition}


\begin{theorem}\label{th_max}
Пусть $\{\psi_\alpha\}_{\alpha\in A}$ --- семейство слабых четностей
в категории диаграмм $\mathfrak K$. Тогда семейство отображений
$\psi_D\colon\V(D)\to\Z_2$, $D\in ob(\mathfrak K)$, где
$\psi_D(v)=\max_{\alpha\in A}(\psi_\alpha)_D(v)$ для каждого
перекрестка $v\in\V(D)$, является слабой четностью. Если все
слабые четности $\psi_\alpha,\ \alpha\in A,$ нетривиальны, то
слабая четность $\psi$ тоже нетривиальна.
\end{theorem}
\begin{proof}
Нам нужно проверить для $\psi$ выполнение свойств $(\Psi0)-(\Psi3)$.

\psz. Пусть $f\colon D\to D'$ --- движение Рейдемейстера и вершина
$v\in\V(D_1)$ не участвует в этом движении. Тогда для любого
$\alpha\in A$ имеет место равенство
$(\psi_\alpha)_{D}(v)=(\psi_\alpha)_{D'}(f_*(v))$. Следовательно,
$$\psi_{D}(v)=\max_{\alpha\in A}(\psi_\alpha)_D(v)=\max_{\alpha\in
A}(\psi_\alpha)_{D'}(f_*(v))=\psi_{D'}(f_*(v)).
$$

\pst. Пусть $f\colon D\to D'$ --- уменьшающее второе движение
Рейдемейстера, в котором участвуют вершины $u,v\in\V(D)$. Тогда
$(\psi_\alpha)_D(u)=(\psi_\alpha)_D(v)$ для любого $\alpha\in A$.
Следовательно,
$$
\psi_{D}(u)=\max_{\alpha\in A}(\psi_\alpha)_D(u)=\max_{\alpha\in
A}(\psi_\alpha)_D(v)=\psi_{D}(v).
$$

\psth. Пусть $f\colon D\to D'$ --- третье движение Рейдемейстера, в
котором участвуют вершины $u,v,w\in\V(D)$. Тогда, очевидно,
$\psi_D(u)=\psi_{D'}(f_*(u))$, $\psi_D(v)=\psi_{D'}(f_*(v))$,
$\psi_D(w)=\psi_{D'}(f_*(w))$. Предположим, что
$\psi_D(u)+\psi_D(v)+\psi_D(w)= 1\in\Z$. Без ограничения общности
можно считать, что $\psi_D(u)=\psi_D(v)=0$ и $\psi_D(w)=1$. Тогда
для любого $\alpha\in A$ имеем
$(\psi_\alpha)_D(u)=(\psi_\alpha)_D(v)=0$. Из свойства \psth для
слабой четности $\psi_\alpha$ следует, что $(\psi_\alpha)_D(w)=0$.
Но тогда $ \psi_{D}(w)=\max_{\alpha\in A}(\psi_\alpha)_D(w)=0, $ и
мы приходим к противоречию.

\pso. Пусть все слабые четности $\psi_\alpha,\ \alpha\in A,$
нетривиальны. Рассмотрим диаграмму $D$ и перекресток $v$ на ней,
который получается первым движением Рейдемейстера. По
теореме~\ref{th_psi1} $(\psi_\alpha)_D(v)=0$ для любого $\alpha\in
A$. Следовательно, $\psi_D(v)=0$, и слабая четность $\psi$
нетривиальна.
\end{proof}

\begin{definition}
Слабую четность $\psi$ из теоремы~\ref{th_max} назовем {\em
максимумом на множестве слабых четностей $\{\psi_\alpha\}_{\alpha\in
A}$} и обозначим как $\max_{\alpha\in A}\psi_\alpha$.
\end{definition}

Отношение порядка на слабых четностях позволяет выделить среди них
наименьший и наибольший элементы. Заметим, что минимальной
слабой четностью является нулевая слабая четность, поскольку
$\mathbb{O}\le\psi$ для любой слабой четности $\psi$. Наибольшим
элементом среди всех слабых четностей является тривиальная
слабая четность $\mathbbm 1$. Рассмотрим теперь только нетривиальные
слабые четности.

\begin{definition}
Нетривиальную слабую четность $\psi_{max}$ в категории диаграмм
$\mathfrak K$ назовем {\em максимальной}, если для любой
нетривиальной слабой четности $\psi$ на $\mathfrak K$ верно
неравенство $\psi\le\psi_{max}$.
\end{definition}

Следующее утверждение очевидно следует из теоремы~\ref{th_max}:
достаточно в качестве семейства слабых четностей
$\{\psi_\alpha\}_{\alpha\in A}$ взять множество всех нетривиальных
слабых четностей в категории $\mathfrak K$.

\begin{corollary}
В любой категории диаграмм $\mathfrak K$ существует единственная
максимальная нетривиальная слабая четность.
\end{corollary}

Далее мы дадим описание максимальной слабой четности для некоторых
категорий диаграмм.

Интерпретация слабых четностей как функториальных отображений дает
возможность определить умножение слабых четностей.
\begin{definition}
Пусть $\psi,\psi'$ --- слабые четности на виртуальных/плоских узлах,
и $\Psi, \Psi'$ --- соответствующие им функториальные отображения.
Тогда слабая четность $\psi\circ\psi'$, соответствующую
функториальному отображению $\Psi\circ\Psi'$, назовем {\em
произведением}
 слабых четностей $\psi,\psi'$.
\end{definition}

\begin{rema}
Несложно проверить, что $\psi\circ\psi'\ge\psi'$ для любых слабых
четностей $\psi,\psi'$. В частности, отсюда следует, что
$\psi\circ\psi_{max}=\psi_{max}$ для любой нетривиальной слабой
четности $\psi$.
\end{rema}

\section{Слабые четности для узлов на фиксированной поверхности}

Пусть $\mathfrak K$ --- категория диаграмм узла (плоского узла)
$\mathcal K$, лежащего на фиксированной двумерной  связной замкнутой
ориентированной поверхности $S$.

\begin{example}[Гомотопическая слабая четность]
Определим {\em гомотопическую слабую четность} $\psi^{hom}$ следующим
образом: для диаграммы $D$ и перекрестка $v\in\V(D)$ положим
$\psi^{hom}_D(v)=0$ в том и только том случае, когда в
гомотопической группе $\pi_1(S,v)$ выполнено соотношение
\begin{equation}\label{hom_pp_eq}
[D_v]=[D]^l
\end{equation}
для некоторого целого числа $l$. Здесь $D_v$
--- одна из половинок узла в вершине $v$ (см. рис.~\ref{pic:halves}). Иными словами,
$\psi^{hom}_D(v)=0$, если гомотопический класс половинки $D_v$ лежит
в подгруппе, порожденной гомотопическим классом узла $\mathcal K$.

\begin {figure}[h]
\centering
\includegraphics[height=2cm]{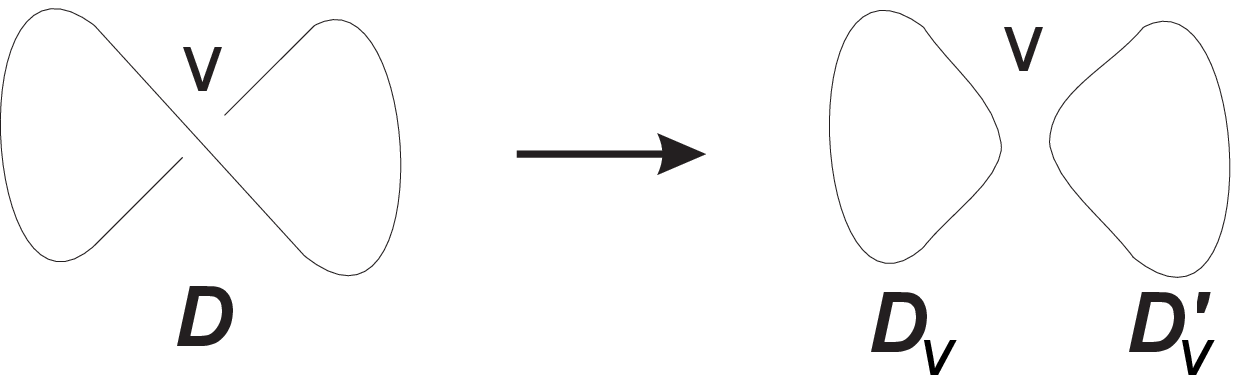}
\caption{Половинки узла в перекрестке}\label{pic:halves}
\end {figure}

Заметим, что классы $[D_v]$ и $[D]$ определены неоднозначно: у узла
есть две половинки $D_v$ и $D'_v$ в вершине $v$, с другой стороны,
гомотопический класс $[\mathcal K]$ равен одному из двух
произведений $[D_v][D'_v]$ или $[D'_v][D_v]$. Таким образом, имеется
четыре варианта равенства~\eqref{hom_pp_eq}. Тем не менее, если
равенство справедливо для какого-либо одного варианта, то оно верно
и для всех остальных (возможно, с другими значениями степени $l$).

Например, пусть $[D_v]=\left([D_v][D'_v]\right)^l$. Тогда
$[D'_v]=[D_v]^{-1}[D_v][D'_v]=\left([D_v][D'_v]\right)^{1-l}$. В
этом случае элементы $[D_v]$ и $[D'_v]$ коммутируют, и
$[D_v][D'_v]=[D'_v][D_v]$.

\begin{prop}
Семейство отображений $\psi^{hom}$ является слабой четностью.
\end{prop}
\begin{proof}
Пусть $D$ --- некоторая диаграмма и $v$ --- ее перекресток.
Соотношение~\eqref{hom_pp_eq} не меняется при гомотопии кривой, образом
которой является диаграмма $D$, если эта гомотопия не затрагивает
перекресток $v$. Следовательно, для $\psi^{hom}$ выполнено свойство
($\Psi0$).

Проверим свойство ($\Psi2$). Пусть перекрестки $u$ и $v$ диаграммы
$D$ участвуют во втором движении Рейдемейстера. Так как двуугольник,
образованный вершинами $u$ и $v$ стягиваем, то мы можем отождествить
гомотопические группы $\pi_1(S,u)$  и $\pi_1(S,v)$. При этом
половинки узла $D_u$ и $D'_u$ отождествляются к половинками $D_v$ и
$D'_v$. Следовательно, соотношение $[D_u]=[D]^l$ эквивалентно
соотношению $[D_v]=[D]^l$. Таким образом, четности перекрестков $u$
и $v$ относительно $\psi^{hom}$ совпадают.

Пусть $u$, $v$ и $w$ --- перекрестки диаграммы $D$, участвующие в
третьем движении Рейдемейстера. Обозначим (в порядке обхода узла)
длинные дуги диаграммы $D$, соединяющие перекрестки $u$, $v$ и $w$,
как $\alpha$, $\beta$ и $\gamma$ (см. рис.~\ref{pic:hom_parity_p3}).
Треугольник, образованный перекрестками, можно стянуть в точку $x\in
S$. При этом $\alpha$, $\beta$ и $\gamma$ отождествляются с
половинками $D_u$, $D_v$ и $D_w$ соответственно, и это соответствие
не зависит от порядка соединения вершин $u$, $v$ и $w$ в диаграмме
$D$. Ниже мы будем использовать обозначения $\alpha$, $\beta$ и
$\gamma$ для гомотопических классов соответствующих кривых. Эти
классы являются элементами группы $\pi_1(S,x)$.

\begin {figure}[h]
\centering
\includegraphics[height=5cm]{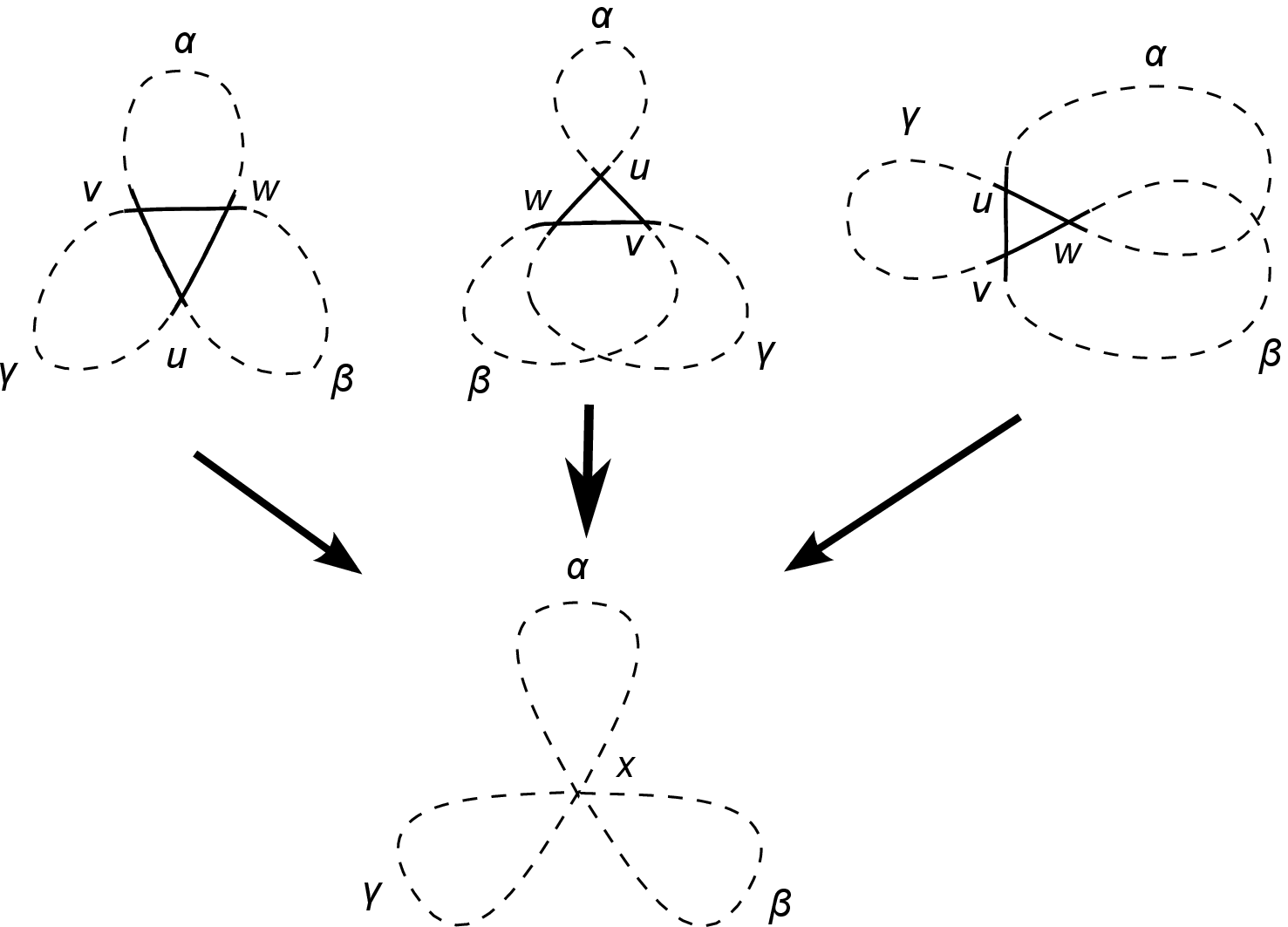}
\caption{Стягивание треугольника в точку при различных способах
соединения его вершин в диаграмме $D$}\label{pic:hom_parity_p3}
\end {figure}

Предположим, что среди перекрестков $u$, $v$ и $w$ один нечетный и
два четных. Не ограничивая общности, можно считать, что нечетным
перекрестком является $w$. Половинки узла в перекрестке $u$ равны
соответственно $\alpha$ и $\beta\gamma$. Так как перекресток $u$
четный, имеет место соотношение $\alpha = (\alpha\beta\gamma)^l$ для
некоторого $l\in\Z$. Кроме того, как было замечено выше, классы
половинок коммутируют, так что
$\alpha\beta\gamma=\beta\gamma\alpha$. Аналогично, для перекрестка
$v$ выполняются соотношения $\beta = (\beta\gamma\alpha)^m, m\in\Z$,
и $\beta\gamma\alpha=\gamma\alpha\beta$. Следовательно, мы имеем
$$
\alpha\beta=(\alpha\beta\gamma)^l(\beta\gamma\alpha)^m=(\gamma\alpha\beta)^{l+m},
$$
откуда $\gamma=(\gamma\alpha\beta)^{1-l-m}$. Таким образом,
перекресток $w$ тоже должен быть четным. Полученное противоречие
показывает, что семейство отображений $\psi^{hom}$ обладает
свойством ($\Psi2$), а значит, является слабой четностью.
\end{proof}
\end{example}

\begin{theorem}
Гомотопическая четность является максимальной нетривиальной слабой четностью.
\end{theorem}

Ключевой для доказательства теоремы является следующая лемма.

\begin{lemma}\label{lem:even_contr}
Пусть $D$ --- диаграмма узла на поверхности $S$ и $v\in\V(D)$ ---
перекресток, такой что $[D_v]=0\in\pi_1(S,v)$ для одной из половинок
узла в перекрестке $v$. Тогда для нетривиальной любой слабой четности $\psi$ имеем
$\psi_D(v)=0$.
\end{lemma}
\begin{proof}
Если узел $\mathcal K$ является плоским, то есть в категории
$\mathfrak K$ не различаются проходы и переходы, то гомотопию
половинки $D_v$ в точку можно представить как последовательность
движений Рейдемейстера. Значит, половинку $D_v$ движениями
Рейдемейстера можно преобразовать в простую петлю в окрестности
перекрестка $v$ (см. рис.~\ref{pic:contraction_flat}). Но тогда
$\psi(v)$ будет равно $0$ по свойству \pso.

\begin {figure}[h]
\centering
\includegraphics[height=2cm]{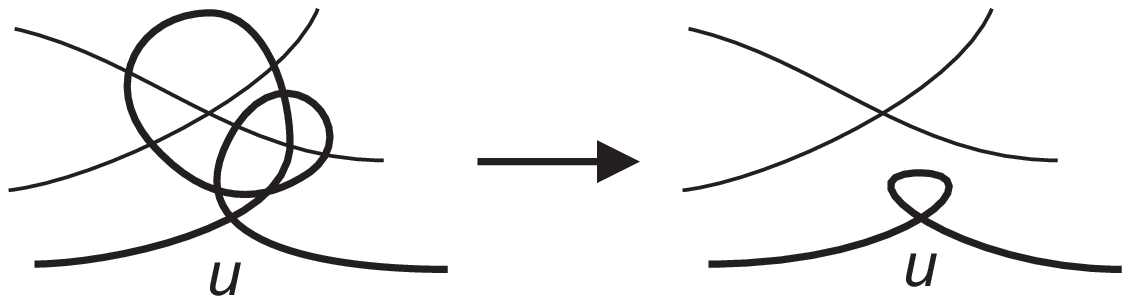}
\caption{Доказательство леммы для плоского узла на
поверхности}\label{pic:contraction_flat}
\end {figure}

Для доказательства леммы в неплоском случае нам потребуются
вспомогательные утверждения.

 \begin{lemma}\label{lem:alt_polygons}
 Для любой нетривиальной слабой четности
 \begin{itemize}
 \item четности вершин в альтернированном двуугольнике совпадают;
 \item число нечетных вершин в альтернированном треугольнике равно $0, 2$ либо $3$.
 \end{itemize}
 \end{lemma}

\begin {figure}[h]
\centering
\includegraphics[height=2cm]{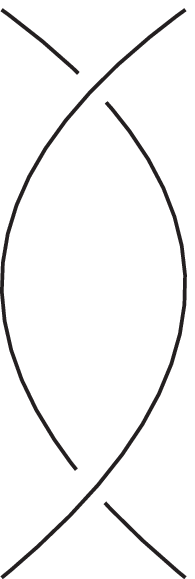}\qquad
\includegraphics[height=2cm]{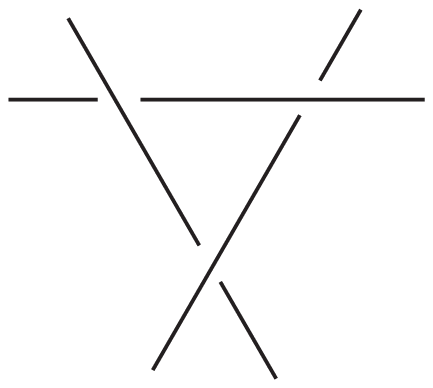}
\caption{Альтернированный двуугольник и
треугольник}\label{pic:alt_poligon}
\end {figure}

\begin{proof}
Пусть перекрестки $u$ и $v$ образуют альтернированный двуугольник.
При помощи первого движения Рейдемейстера образуем на дуге,
соединяющей $u$ и $v$, перекресток $w$, чтобы к перекресткам $u,v$ и
$w$  можно было применить третье движение Рейдемейстера (см.
рис.~\ref{pic:alt_bigon}). Перекресток $w$, согласно свойству
($\Psi1$), является четным, так что четности перекрестков $u$ и $v$
должны совпадать, чтобы выполнялось свойство ($\Psi3$).
\begin {figure}[h]
\centering
\includegraphics[height=1.5cm]{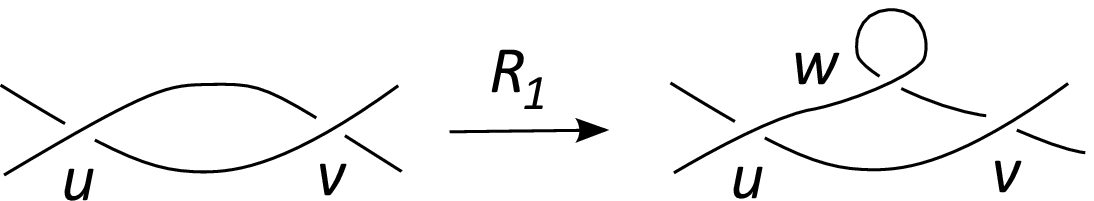}
\caption{Доказательство леммы для альтернированного
двуугольника}\label{pic:alt_bigon}
\end {figure}

Пусть перекрестки $u,v$ и $w$ образуют альтернированный треугольник.
При помощи второго движения Рейдемейстера добавим перекрестки $t$ и
$t'$ (см. рис.~\ref{pic:alt_trigon}). Согласно ($\Psi3$), количество
нечетных перекрестков среди $v,w$ и $t$ может быть равно $0$, $2$
или $3$. Но четности перекрестков $t$ и $t'$ совпадают, согласно
($\Psi2$), а совпадение четностей перекрестков $t'$ и $u$,
образующих альтернированный двуугольник, следует из доказанного
выше.

\begin {figure}[h]
\centering
\includegraphics[height=2.5cm]{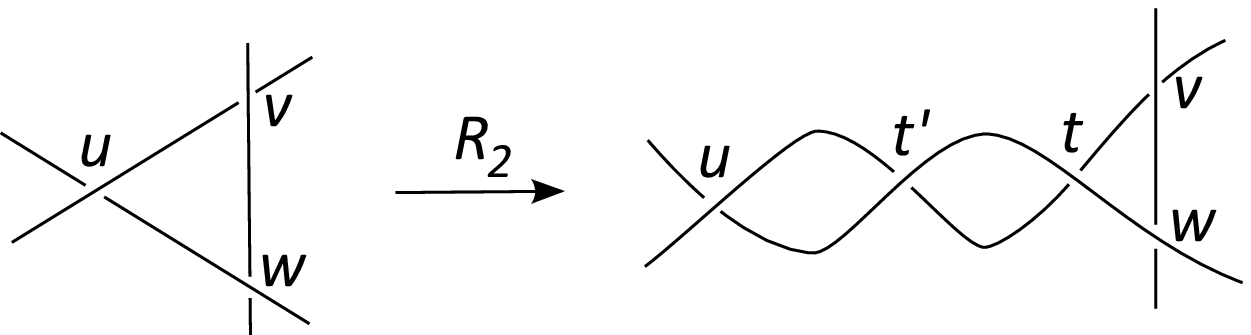}
\caption{Доказательство леммы для альтернированного
треугольника}\label{pic:alt_trigon}
\end {figure}

Таким образом, лемма~\ref{lem:alt_polygons} доказана.
\end{proof}

 \begin{lemma}\label{lem:crossed_polygons}
 Для любой нетривиальной слабой четности
 \begin{itemize}
 \item четности вершин $u,v$ в пересеченном двуугольнике совпадают;
 \item вершина $u$ в пересеченном одноугольнике является четной
 \end{itemize}
при любом выборе проходов и переходов в перекрестках (см.
рис.~\ref{pic:transversed_poligons}).
\end{lemma}

\begin {figure}[h]
\centering
\includegraphics[height=2cm]{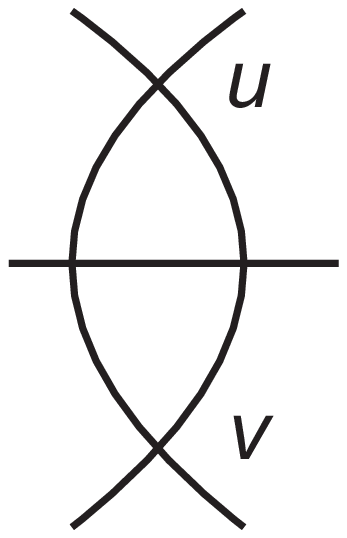}\qquad
\includegraphics[height=2cm]{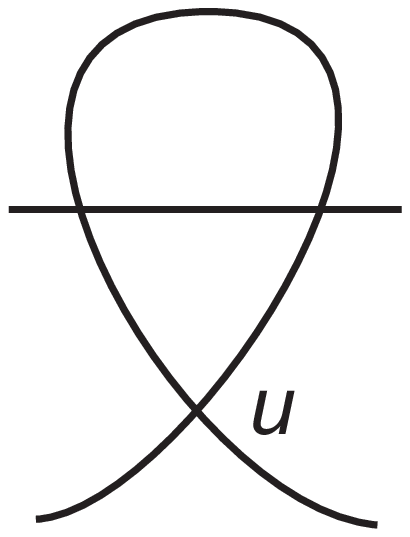}
\caption{Пересеченный двуугольник и
одноугольник}\label{pic:transversed_poligons}
\end {figure}

\begin{proof}
Рассмотрим пересеченный двуугольник с вершинами $u,v$. В случае,
когда пересекающую дугу можно вытащить из двуугольника при помощи
третьего движения, совпадение четностей вершин $u$ и $v$ следует из
свойства ($\Psi0$), а также свойства ($\Psi2$) либо
леммы~\ref{lem:alt_polygons}.

Предположим теперь, что пересекающую дугу из  двуугольника удалить
третьим движением нельзя (см. рис.~\ref{pic:crossed_bigon}). Добавим
при помощи второго движения Рейдемейстера перекрестки $w$ и $w'$.
Проведенные выше рассуждения показывают, что четности перекрестков
$u$ и $w$ совпадают. С другой стороны, четности перекрестков $w$ и
$w'$ совпадают, согласно свойству ($\Psi2$), а четности перекрестков
$w'$ и $v$ совпадают по лемме~\ref{lem:alt_polygons}. Таким образом,
четности перекрестков $u$ и $v$ равны.

\begin {figure}[h]
\centering
\includegraphics[height=2.5cm]{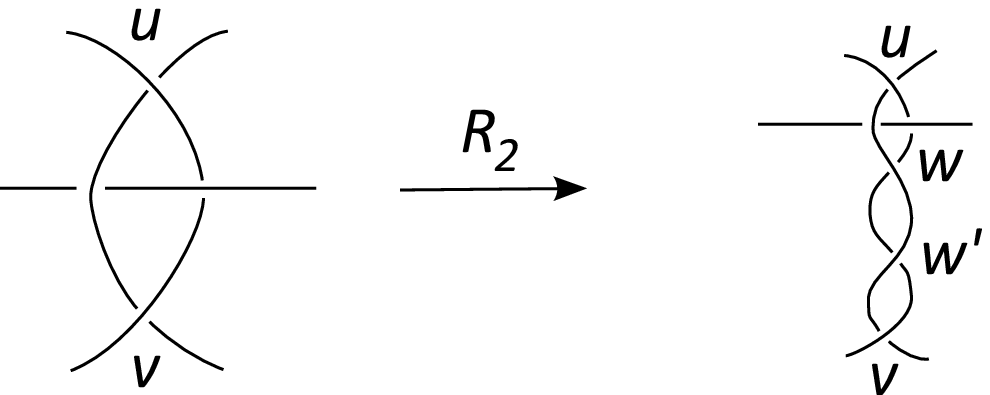}
\caption{Доказательство леммы для пересеченного
двуугольника}\label{pic:crossed_bigon}
\end {figure}

Рассмотрим пересеченный одноугольник с вершиной $u$ (см.
рис.~\ref{pic:transversed_poligons} справа). При помощи первого
движения Рейдемейстера мы можем создать вверху петли перекресток
$v$, который будет четным, согласно свойству ($\Psi1$). Перекрестки
$u$ и $v$ будут образовывать пересеченный двуугольник, поэтому их
четности будут одинаковы. Следовательно, перекресток $u$ четный.
Лемма~\ref{lem:crossed_polygons} доказана.
\end{proof}

Следующая лемма является аналогом леммы 4.7 из работы~\cite{IMN}.
\begin{lemma}\label{lem:parity_polygon}
Пусть диаграмма $D$ ограничивает на поверхности многоугольник,
вершинами которого являются перекрестки $v_0,v_1,\dots, v_n$.
Предположим, что перекрестки $v_1,\dots, v_n$ являются четными
относительно некоторой слабой четности $\psi$. Тогда перекресток
$v_0$  тоже четный относительно $\psi$.
\end{lemma}

\begin{proof}
Доказательство ведется по индукции. При $n=0,1,2$ утверждение
следует из свойств ($\Psi1$), ($\Psi2$) и ($\Psi3$) соответственно
(а также леммы~\ref{lem:alt_polygons}).

Пусть утверждение леммы верно при $n\le k$. Рассмотрим случай
$n=k+1$. Пусть имеется $(k+2)$-угольник $v_0,v_1,\dots, v_{k+1}$, в
котором перекрестки $v_1,\dots, v_{k+1}$  четные. Применим
увеличивающее второе движение Рейдемейстера к ребрам $v_0v_1$ и
$v_kv_{k+1}$, создав два новых перекрестка $w$ и $w'$. На месте
$(k+2)$-угольника появятся треугольник $v_0wv_{k+1}$, двуугольник
$ww'$  и $(k+1)$-угольник $w'v_1\dots v_k$ (см. рис.~\ref{pic:even_polygon}).
\begin {figure}[h]
\centering
\includegraphics[height=2.5cm]{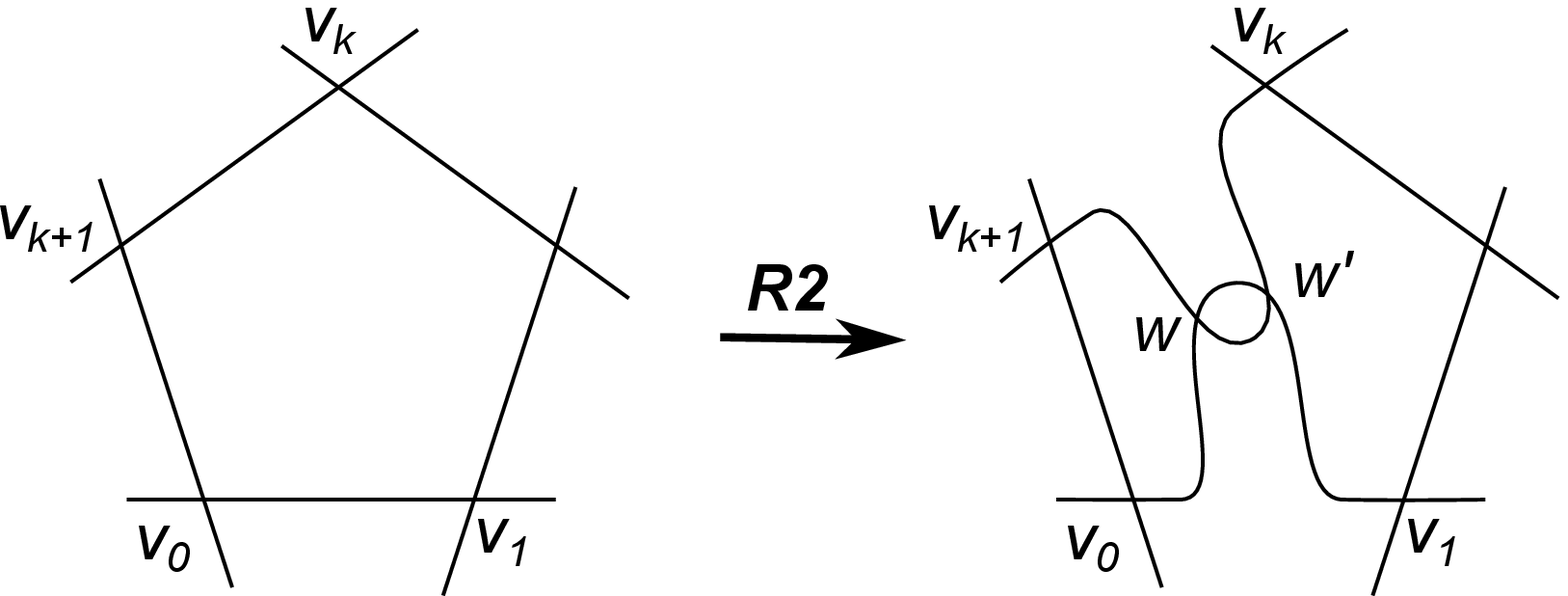}
\caption{Доказательство леммы~\ref{lem:parity_polygon}}\label{pic:even_polygon}
\end {figure}
По предположению индукции
перекресток $w'$ будет четным. Согласно свойству ($\Psi2$) четным
окажется и перекресток $w$. Поскольку $v_{k+1}$ тоже является
четной, то оставшаяся вершина $v_0$ в треугольнике $v_0wv_{k+1}$
будет четной по свойству ($\Psi3$) либо
лемме~\ref{lem:alt_polygons}.

Лемма~\ref{lem:parity_polygon} доказана.
\end{proof}

Завершим доказательство леммы~\ref{lem:even_contr}.

Рассмотрим перекресток $u$ диаграммы $D$, половинка которого $D_u$
стягивается в точку (см. рис.~\ref{pic:contracting0}).
\begin {figure}[h]
\centering
\includegraphics[height=2.5cm]{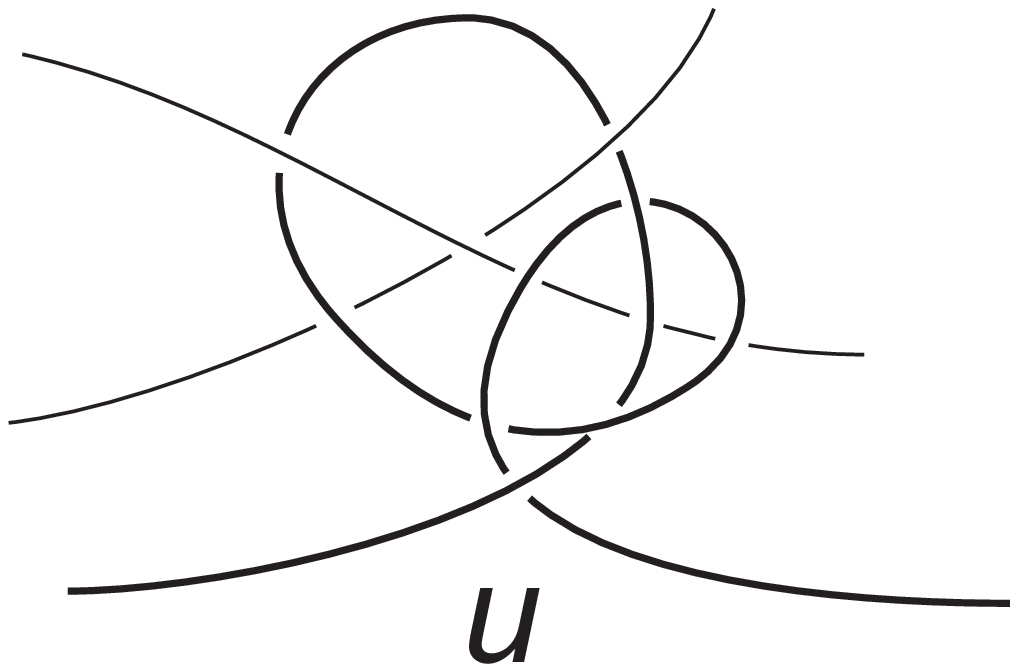}
\caption{Гомотопически тривиальная
половинка}\label{pic:contracting0}
\end {figure}
Процесс стягивания можно представить как проекцию некоторой
гомотопии узла в утолщении поверхности $S\times I$. Данная гомотопия
представляется как последовательность изотопий в $S\times I$ и
моментов самопересечения узла, которые при проекции выглядят как
операции замены прохода на переход в некотором перекрестке (см.
рис.~\ref{pic:cross_change}).
\begin {figure}[h]
\centering
\includegraphics[height=1.5cm]{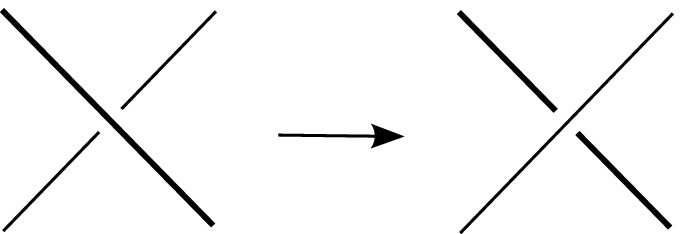}
\caption{Замена прохода на переход в
перекрестке}\label{pic:cross_change}
\end {figure}
Мы модифицируем процесс гомотопии добавляя в момент замены
перекрестка небольшую петлю, чтобы избежать самопересечений (см.
рис.~\ref{pic:cross_change1}).
\begin {figure}[h]
\centering
\includegraphics[height=2cm]{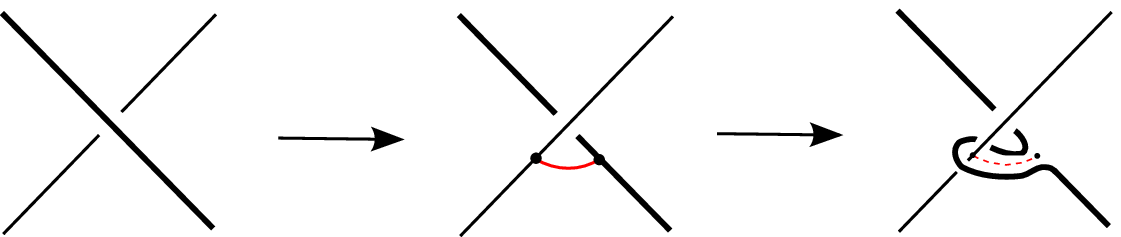}
\caption{Добавление петли вместо замены структуры
перекрестка}\label{pic:cross_change1}
\end {figure}
В результате мы получим изотопию узла, проекция которой переводит
половинку $D_u$ в маленькую петлю с добавленными "отростками" (см.
рис.~\ref{pic:contracting1}).
\begin {figure}[h]
\centering
\includegraphics[height=5cm]{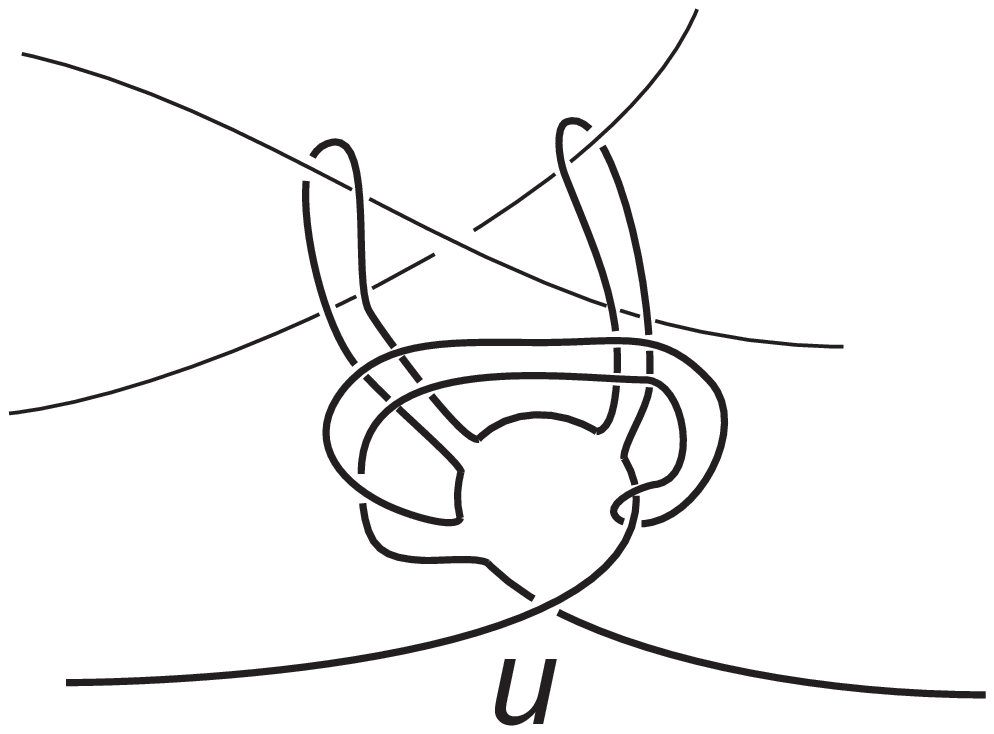}
\caption{Стянутая половинка с отростками}\label{pic:contracting1}
\end {figure}
Применим второе движение Рейдемейстера необходимое число раз, чтобы
образовался маленький многоугольник, содержащий перекресток $u$, все
остальные вершины которого соответствуют "отросткам" (см.
рис.~\ref{pic:contracting2}).
\begin {figure}[!h]
\centering
\includegraphics[height=5cm]{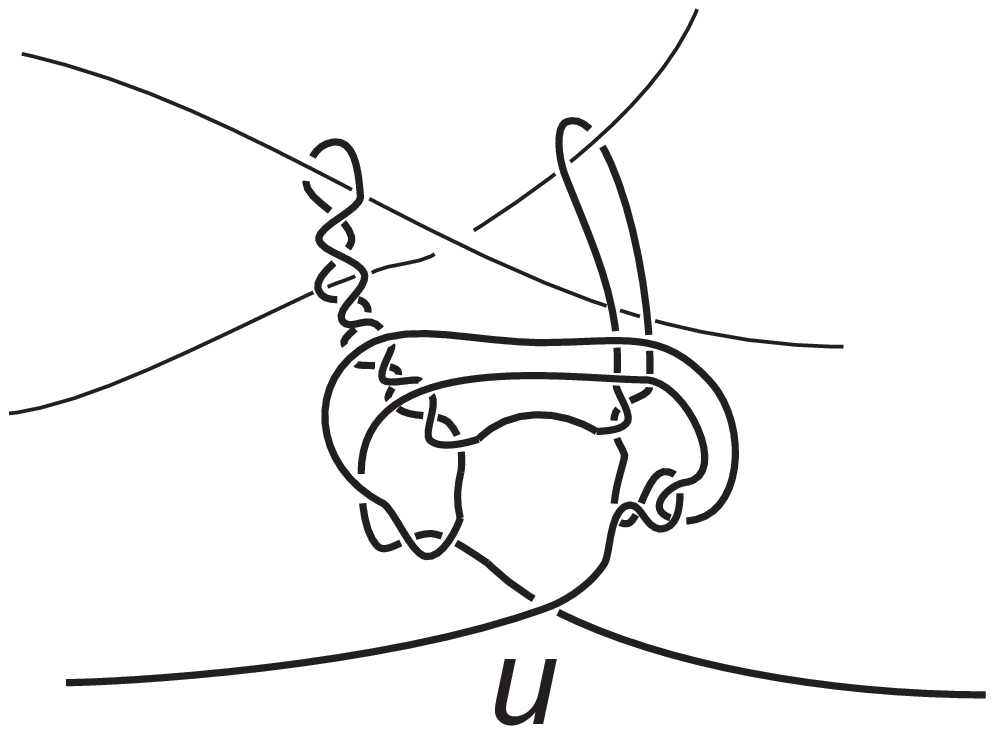}
\caption{Добавление перекрестков в конце отростков}\label{pic:contracting2}
\end {figure}
Заметим, что все эти вершины будут четными, поскольку каждый
отросток можно разбить вторыми движениями на (возможно,
пересеченные) двуугольники (см. левый верхний отросток на рис.~\ref{pic:contracting2}),
так что все перекрестки отростка будут иметь одинаковую четность,
согласно свойству ($\Psi2$) и
леммам~\ref{lem:alt_polygons},\ref{lem:crossed_polygons}. Эта
четность будет совпадать с четностью первого перекрестка отростка,
который является вершиной пересеченного одноугольника, а значит,
четен. Таким образом, перекресток $u$ является вершиной
многоугольника, остальные вершины которого четны. Следовательно, по
лемме~\ref{lem:parity_polygon}, перекресток $u$ четный.

Лемма~\ref{lem:even_contr} доказана.
\end{proof}

\begin{proof}[Доказательство теоремы~\ref{th_max}]
Нам достаточно показать, что для любой слабой четности $\psi$ на
категории диаграмм узла в поверхности $S$
соотношение~\eqref{hom_pp_eq} влечет четность соответствующего
перекрестка.

Обозначим через $L$ множество целых чисел, обладающих следующим
свойством: $l\in L$,  если для любой диаграммы $D$, перекрестка
$v\in\V(D)$, такого что $[D_v]=[D]^l\in\pi_1(S,v)$ для некоторой
половинки $D_v$ узла в $v$, и любой слабой четности $\psi$ имеет
место соотношение $\psi_D(v)=0$.

Согласно лемме~\ref{lem:even_contr} выполнено включение $0\in L$.

Очевидно, что из $l\in L$ следует $1-l\in L$ (это соответствует
переходу от одной половинке к другой). Таким образом, $1=1-0\in L$.

Покажем теперь, что из $l\in L$ следует $l-1\in L$. Пусть $l\in L$ и
для перекрестка $v$ выполнено соотношение $[D_v]=[D]^{l-1}$. Завяжем
рядом с перекрестком $v$ петлю с гомотопическим типом $[D_u]=[D]$
(здесь мы тождествляем фундаментальные группы с базовыми точками,
принадлежащими стягиваемой окрестности точки $v$).  Четность
перекрестка $v$ при этом не изменится согласно свойству ($\Psi0$).
Применим второе движение Рейдемейстера, образовав перекресток $w$.
Гомотопический тип половинки узла в перекрестке $w$ равен
$$
[D_w]=[D_v][D_u]=[D]^{l-1}[D]=[D]^l.
$$
Так как $l\in L$, перекресток $w$ оказывается четным. Поскольку
$1\in L$, перекресток $u$ тоже четен. Но перекрестки $u,v$ и $w$
образуют треугольник. Значит, согласно свойству ($\Psi3$) либо
лемме~\ref{lem:alt_polygons}, перекресток $v$ будет четным.

\begin {figure}[h]
\centering
\includegraphics[height=1.5cm]{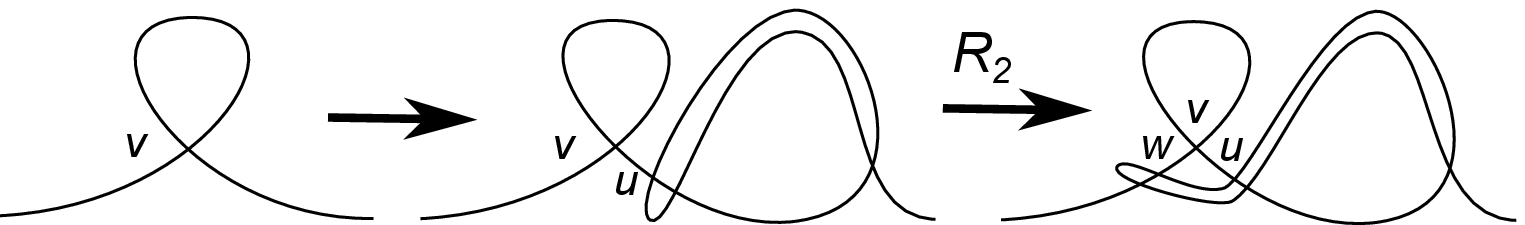}
\caption{Доказательство перехода $l\to l-1$}\label{pic:hom_max}
\end {figure}

Таким образом, все целые числа, не превосходящие $1$ лежат в $L$.
Следовательно, $\N\in L$, так как для любого $n\in\N$ выполнено $n =
1-(1-n)$ и $1-n<1$. Тогда $L=\Z$.

Итак, любой перекресток, четный относительно гомотопической
четности, является четным для любой нетривиальной слабой четности.
Значит, гомотопическая четность $\psi^{hom}$ является максимальной.
\end{proof}

\begin{corollary}\label{cor:classical_weak_parities}
На классических узлах любая слабая четность является нулевой либо
тривиальной.
\end{corollary}
\begin{proof}
Классические узлы суть узлы на сфере. Так как фундаментальная группа
сферы тривиальна, то все перекрестки являются четными относительно
гомотопической четности, то есть $\psi^{hom}= \mathbb{O}$. Так как
гомотопическая четность максимальна, любая нетривиальная слабая
четность совпадает с нулевой.
\end{proof}

Таким образом, как и в случае четностей (см.~\cite[Corollary
4.2]{IMN}), теория слабых четностей становится тривиальной на
классических узлах.

Переформулировка следствия~\ref{cor:classical_weak_parities} на язык
функториальных отображений приводит к следующему результату.

\begin{corollary}
Ограничение любого нетривиального функториального отображения на классические узлы
является тождественным отображением. В частности, нетривиальное функториальное
отображение виртуальных узлов, принимающее значение в классических
диаграммах узлов, есть проекция на классические узлы.
\end{corollary}

Опишем теперь все слабые четности для узлов на заданной двумерной
поверхности. При этом мы наложим на диаграммы узла некоторое
необременительное ограничение -- - потребуем, чтобы диаграммы
проходили через фиксированную точку на поверхности.

Пусть $\mathcal K$ --- узел в утолщении связной ориентированной
замкнутой поверхности $S$, а $z$ --- некоторая точка на поверхности
$S$. Рассмотрим категорию диаграмм $\mathfrak K_z$, объектами
которой являются диаграммы узла $\mathcal K$, содержащие точку $z$,
причем $z$ не является перекрестком. Морфизмы категории --- это
изотопии и движения Рейдемейстра, оставляющие точку $z$ неподвижной,
а также их композиции.

\begin {figure}[h]
\centering
\includegraphics[height=2cm]{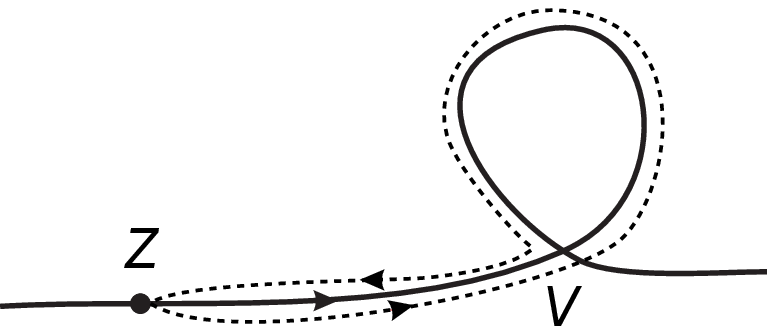}
\caption{Приведенная половинка диаграммы в перекрестке
$v$}\label{pic:punctured_half}
\end {figure}

 \begin{theorem}
Имеется естественная биекция между слабыми четностями на категории
диаграмм $\mathfrak K_z$ и подгруппами фундаментальной группы
$H\subset \pi_1(S,z)$, такими, что $[\mathcal K]\in H$.

При заданной подгруппе $H$ соответствующая слабая четность $\psi^H$
определяется следующим образом: перекресток $v$ диаграммы $D$
является четным, т.е.$\psi^H_D(v)=0$, тогда и только тогда, когда
$[\hat D_{v}]\in H$, где $\hat
 D_{v}$ --- это приведенная к $z$ половинка диаграммы в перекрестке $v$ (см. рис.~\ref{pic:punctured_half}).
 \end{theorem}

\begin{rema}
Определение приведенной половинки зависит от ориентации узла,
однако, как мы покажем ниже при доказательстве теоремы, на
выполнение условия четности перекрестка выбор ориентации узла
влияния не оказывает.
\end{rema}

\begin{proof}
Пусть $H$ --- подгруппа в $\pi_1(S,z)$, содержащая элемент
$[\mathcal K]$. Покажем, что семейство отображений $\psi^H =
\{\psi^H_D\}$,
$\psi^H_D\colon\V(D)\to\Z_2$, $D\in ob(\mathfrak K_z)$, определяемое условием
 $$
\psi^H_D(v)=\left\{
              \begin{array}{cl}
                0, & [\hat D_{v}]\in H; \\
                1, & \hbox{иначе},
              \end{array}
            \right.
 $$
является слабой четностью.

Покажем сначала, что четность перекрестка не зависит от того,
гомотопический класс какой половинки проверяется на принадлежности
подгруппе $H$, а также от ориентации узла.

Рассмотрим перекресток $v$ и обозначим дугу узла, соединяющую (в
направлении обохода узла) точки $z$ и $v$, как $\gamma_1$;
(неприведенную) половинку узла в перекрестке $v$, не содержащую $z$,
как $\alpha$; и дугу, соединяющую (по направлению обхода узла) $v$ и
$z$, как $\gamma_2$ (см. рис.~\ref{pic:based_diagram}).
\begin {figure}[h]
\centering
\includegraphics[height=3cm]{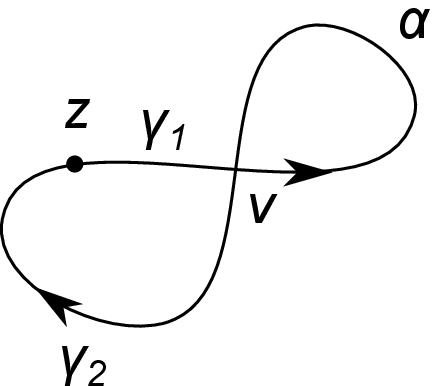}
\caption{Перекресток в диаграмме с отмеченной точкой $z$}\label{pic:based_diagram}
\end {figure}
Тогда гомотопический класс одной приведенной половинки в перекрестке
$v$ равен $[\hat D_v]=[\gamma_1\alpha\gamma_1^{-1}]$, а другой ---
равен $[\hat
D'_v]=[\gamma_1\alpha\gamma_2\gamma_1\alpha^{-1}\gamma_1^{-1}]=[D]
[\hat D_v^{-1}]$, где $[D]=[\gamma_1\alpha\gamma_2]$ ---
гомотопический класс диаграммы (и узла). Так как $[D]\in H$,
половинки $[\hat D_v]$ и $[\hat D'_v]$ принадлежат либо не
принадлежат подгруппе $H$ одновременно.

Приведенная половинка диаграммы в перекрестке при противоположной
ориентации узла равна
${(-\hat D)}_v=\gamma_2^{-1}\alpha^{-1}\gamma_2$. Тогда имеет место
равенство $[{(-\hat D)}_v]=[D]^{-1}[\hat D_v]^{-1}[D]$. Так как
$[D]\in H$, условие $[\hat D_v]\in H$ эквивалентно условию
$[{(-\hat D)}_v]\in H$.

Проверим теперь выполнение условий слабой четности.

Свойство ($\Psi0$) следует из того, что гомотопический класс
половинки диаграммы в перекрестке не меняется при изотопии и
движениях Рейдемейстера, при которых перекресток сохраняется, однако
следует проверить ситуацию, когда перекресток диаграммы $v$
переходит через точку $z$ (см. рис.~\ref{pic:over_based_point}).
\begin {figure}[h]
\centering
\includegraphics[height=3cm]{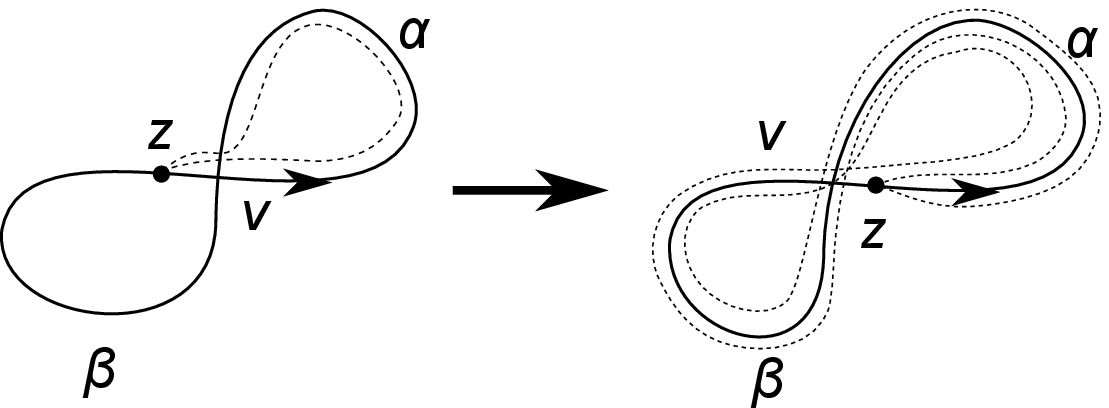}
\caption{Переход перекрестка через отмеченную точку
$z$}\label{pic:over_based_point}
\end {figure}
Обозначим через $\alpha$ и $\beta$ (неприведенные) половинки
диаграммы в перекрестке $v$. Короткую дугу, соединяющую $v$ и $z$,
можно стянуть в точку; при этом половинки диаграммы до перехода
перекрестка через $z$ можно отождествить с соответствующими
половинками после перехода и рассматривать как петли с началом в
точке $z$. Приведенная половинка $\hat D_v$ в начальной диаграмме
равна $\alpha$, после перехода ей соответствует приведенная
половинка $\hat D'_v=\alpha\beta\alpha\beta^{-1}\alpha^{-1}$ (см.
рис.~\ref{pic:over_based_point}). Тогда $[\hat D'_v]=[D][\hat
D_v][D]^{-1}$. Так как $[D]\in H$, четность перекрестка при переходе
через $z$ не меняется.

Проверим свойство ($\Psi2$). Пусть вершины $u$ и $v$ принимают
участие во втором движении Рейдемейстера (см.
рис.~\ref{pic:h_parity_p2}). Тогда $\hat D_u =
\gamma_1\delta_1\alpha\delta_2\gamma_1^{-1}$ и
$\hat D_v=\gamma_1\delta_1\alpha\delta_1^{-1}\gamma_1^{-1}$. Так как
малые дуги $\delta_2$ и $\delta_1^{-1}$ гомотопны как кривые с
закрепленными концами, имеет место равенство $[\hat D_u]=[\hat
D_v]$,  следовательно, перекрестки $u$ и $v$ имеют одинаковую
четность относительно $\psi^H$.
\begin {figure}[h]
\centering
\includegraphics[height=3cm]{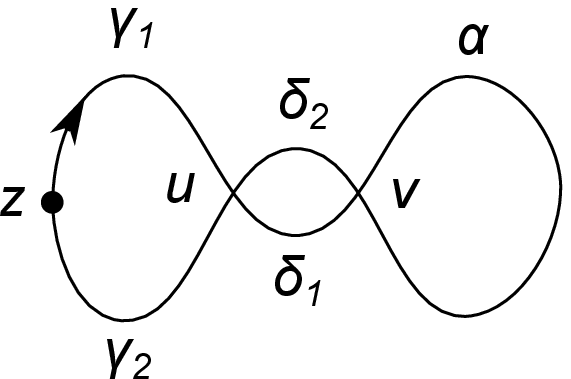}
\caption{Доказательство свойства ($\Psi2$)}\label{pic:h_parity_p2}
\end {figure}

Пусть перекрестки $u,v$ и $w$ участвуют в третьем движении (см.
рис.~\ref{pic:h_parity_p3}).
\begin {figure}[h]
\centering
\includegraphics[height=3cm]{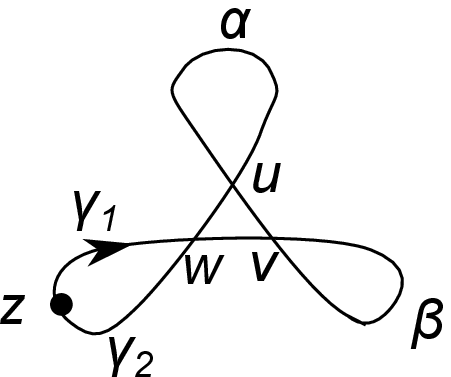}
\caption{Доказательство свойства ($\Psi3$)}\label{pic:h_parity_p3}
\end {figure}
Тогда их приведенные половинки равны $\hat
D_u=\gamma_1\alpha\gamma_1^{-1}$, $\hat
D_v=\gamma_1\alpha\beta\alpha^{-1}\gamma_1^{-1}$ и $\hat
D_w=\gamma_1\alpha\beta\gamma_2\gamma_1\beta^{-1}\alpha^{-1}\gamma_1^{-1}$.
Следовательно, $[\hat D_w][\hat D_v][\hat
D_u]=[\gamma_1\alpha\beta\gamma_2]=[D]\in H$. Отсюда вытекает, что
количество половинок перекрестков, гомотопический класс которых
принадлежит подгруппе $H$, может быть равным $0,1$ или $3$. Иными
словами, количество нечетных перекрестков среди $u,v$ и $w$ не может
быть равно $1$, что доказывает свойство ($\Psi3$). Проверка свойства
($\Psi3$) для других конфигураций перекрестков производится
аналогично.

Таким образом, мы показали, что $\psi^H$ является слабой четностью.
Проверим теперь, что любая слабая четность на категории диаграмм
$\mathfrak K_z$ равна $\psi^H$ для некоторой подгруппы $H$.

Пусть $\psi$  --- некоторая слабая четность на $\mathfrak K_z$.

Предположим, что диаграмма $D$ содержит перекрестки $u$ и $v$, такие
что $[\hat D_u]=[\hat D_v]\in\pi_1(S,z)$. Покажем, что тогда
$\psi_D(u)=\psi_D(v)$. Для этого притянем перекрестки $u$ и $v$ к
отмеченной точке (см. рис.~\ref{pic:base_reduction}), так чтобы они
попали в стягиваемую окрестность точки $z$.
\begin {figure}[h]
\centering
\includegraphics[height=1.5cm]{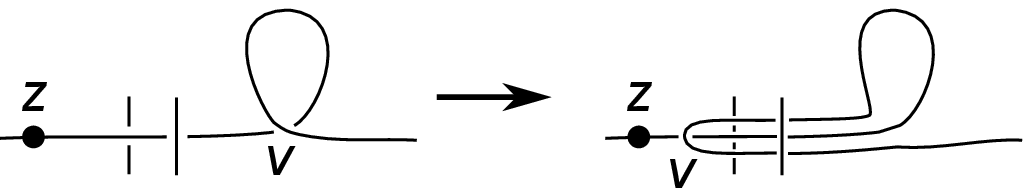}
\caption{Притягивание перекрестка к отмеченной
точке}\label{pic:base_reduction}
\end {figure}
Заметим, что притягивание перекрестка не меняет гомотопический тип
его приведенной половинки.

Обозначим в получившейся диаграмме через $\alpha$ половинку
диаграммы в перекрестке $u$, не содержащую $z$,  через $\beta$ --- длинную дугу
диаграммы, соединяющую точки $u$ и $v$, и через $\gamma$ --- длинную
дугу, соединяющую $v$ и $z$ (см. рис.~\ref{pic:h_parity_hom_inv}).
\begin {figure}[h]
\centering
\includegraphics[height=5cm]{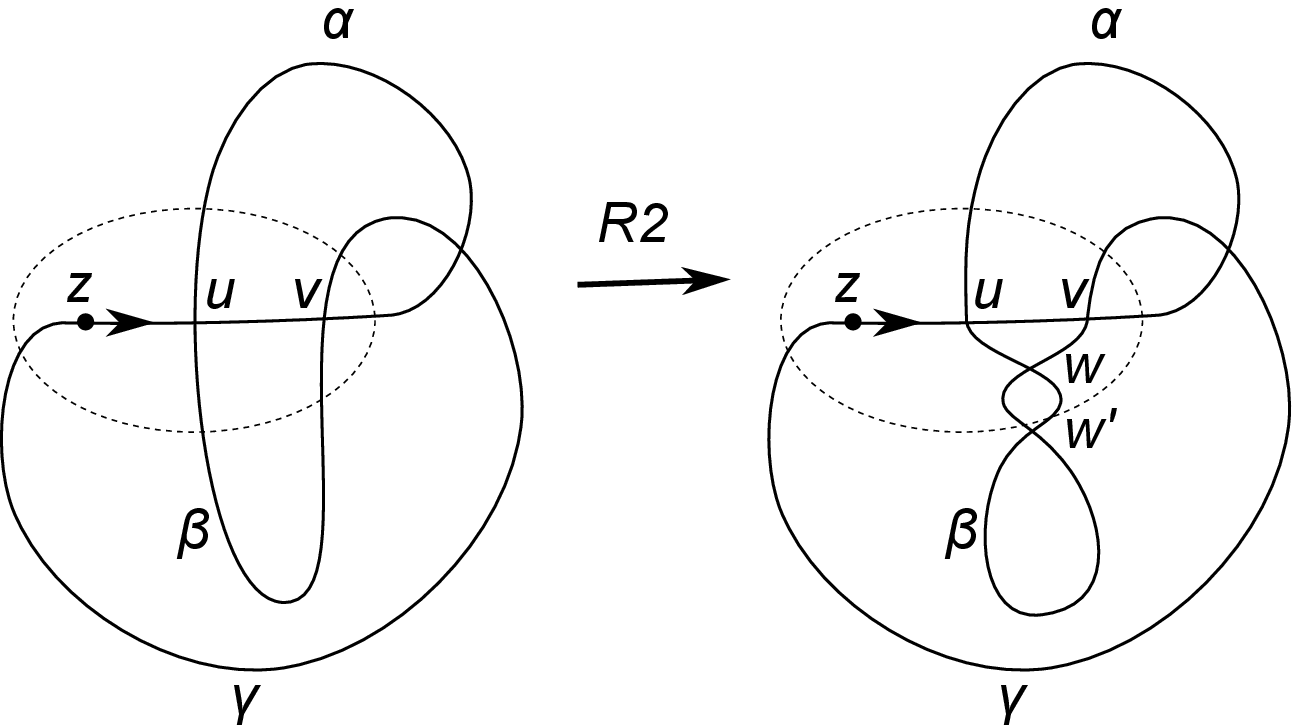}
\caption{Диаграмма с притянутыми к отмеченной точке перекрестками
$u$ и $v$. Пунктиром обозначена стягиваемая окрестность точки $z$}\label{pic:h_parity_hom_inv}
\end {figure}
Отождествляя между собой точки стягиваемой окрестности, мы будет
игнорировать дуги, целиком лежащие в этой окрестности, и
рассматривать кривые $\alpha,\beta,\gamma$ как замкнутые петли с
началом в $z$. Тогда приведенная половинка в перекрестке $u$ равна
$\hat D_u=\alpha$, а половинка в перекрестке $v$ равна $\hat
D_v=\alpha\beta$. Так как по предположению $[\hat D_u]=[\hat D_v]$,
то $[\beta]=1\in\pi_1(S,z)$. Применим к диаграмме второе движение
Рейдемейстера, добавляющее перекрестки $w$ и $w'$. Одна из
(неприведенных) половинок перекрестка $w$ равна $\beta$, а значит,
стягиваема. По лемме~\ref{lem:even_contr} имеем $\psi_D(w)=0$.
Тогда, так как перекрестки $u,v$ и $w$ образуют треугольник, из
свойства ($\Psi3$) либо леммы~\ref{lem:alt_polygons} следует, что
перекрестки $u$ и $v$ имеют одинаковую четность.

Пусть теперь имеются перекрестки $u\in\V(D)$ и $v\in\V(D')$,
$D,D'\in ob(\mathfrak K_z)$, такие что $[\hat D_u]=[\hat
D'_v]\in\pi_1(S,z)$. Так как диаграммы $D$ и $D'$ соответствуют
одному узлу $\mathcal K$, имеется последовательность движений,
преобразующих $D$ в $D'$. При этих преобразованиях перекресток $u$
может исчезнуть, но мы можем при помощи вторых движений
Рейдемейстера добавить к диаграмме  $D$ перекресток $w$, такой что
$[\hat D_w]=[\hat D_u]$, причем этот перекресток будет сохраняться
при переходе от $D$ к $D'$ (после некоторой модификации
промежуточных диаграмм, состоящей в добавлении к ним перекрестка $w$
при помощи вторых и, возможно, третьих движений Рейдемейстера).
Гомотопический тип половинки и четность перекрестка $w$ при
движениях сохраняются, поэтому
$\psi_D(u)=\psi_D(w)=\psi_{D'}(w)=\psi_{D'}(v)$, где первое и
последнее равенства следуют из рассуждений, проведенных ранее.

Определим множество $H\subset\pi_1(S,z)$, состоящее из
гомотопических классов приведенных половинок в перекрестках, четных
относительно $\psi$:
 $$
H=\{[\hat D_v]\,|\, D\in ob(\mathfrak K_z), v\in\V(D) \mbox{ такие,
что } \psi_D(v)=0\}.
 $$
Тогда проведенные выше рассуждения показывают, что слабая четность
$\psi$ задается по правилу
 $$
\psi_D(v)=\left\{
              \begin{array}{cl}
                0, & [\hat D_{v}]\in H; \\
                1, & \hbox{иначе},
              \end{array}
            \right.
 $$
Из свойства ($\Psi1$) следует, что $[\mathcal K]\in H$. Нам остается
показать, что $H$ является подгруппой в $\pi_1(S,z)$.

Пусть $h_1,h_2\in H$. На произвольной диаграмме $D$ при помощи
вторых движений Рейдемейстера образуем перекрестки $u$ и $v$, такие
что $[\hat D_u]=h_1$, $[\hat D_v]=h_2$. Тогда
$\psi_D(u)=\psi_D(v)=0$. Притянем перекрестки $u$ и $v$ к отмеченной
точке $z$ и при помощи второго движения добавим перекрестки $w$ и
$w'$ (см. рис.~\ref{pic:h_parity_hom_inv}). Тогда имеем $\hat
D_u=\alpha$, $\hat D_v=\alpha\beta$, $\hat
D_w=\alpha\beta\alpha^{-1}$ (здесь мы опять не обращаем внимания на
дуги, целиком лежащие внутри стягиваемой окрестности точки $z$).
Следовательно, $[\hat D_w]=[\hat D_v][\hat D_u]^{-1}=h_2h_1^{-1}$. В
силу свойства ($\Psi3$) (либо леммы~\ref{lem:alt_polygons})
перекресток $w$ является четным относительно $\psi$, а значит,
$h_2h_1^{-1}\in H$. Поскольку это включение имеет место для любых
$h_1,h_2\in H$, множество $H$ является подгруппой. Таким образом,
$\psi=\psi^H$.
\end{proof}

\begin{rema}
Слабую четность, определенную на семействе диаграмм узла $\mathcal
K$,  проходящих через точку $z$, можно распространить на множество
всех диаграмм узла $\mathfrak K$ следующим образом. Фиксируем
диаграмму $D_0\in\mathfrak K_z$. Пусть $D$ --- произвольная
диаграмма узла $\mathcal K$. Так как диаграммы $D_0$ и $D$
соответствуют одному узлу, имеется изотопия в утолщении поверхности,
проекция которой на поверхность переводит диаграмму $D_0$ в $D$. При
изотопии точка $z$, рассматриваемая как точка на диаграмме, опишет
путь $\gamma$ с началом в $z$ и концом в точке $w\in D$.

Пусть слабая четность $\psi$ в $\mathfrak K_z$ определяется
подгруппой $H\subset\pi_1(S,z)$. Путь $\gamma$ определяет изоморфизм
фундаментальных групп $\gamma_*\colon \pi_1(S,z)\to\pi_1(S,w)$.
Значение слабой четности $\psi$ в перекрестке $v$ диаграммы $D$
может быть определено по правилу: $\psi_D(v)=0$ тогда и только
тогда, когда $[\hat D_v]\in\gamma_*(H)\subset\pi_1(S,w)$, где $\hat
D_v$
--- половинка узла в перекрестке $v$, приведенная к точке $w$.
\end{rema}

Пусть $\mathcal K$ --- узел на поверхности $S$ и $D$ --- диаграмма
узла $\mathcal K$, содержащая точку $z$. Пусть $H$ --- подгруппа в
$\pi_1(S,z)$, которая содержит класс $[\mathcal K]$.  Подгруппе $H$
соответствует накрытие $p_H\colon \widetilde S\to S$ поверхности
$S$. При этом поднятие $\widetilde D$ диаграммы $D$, содержащее
отмеченную точку $\widetilde z$ поверхности $\widetilde S$
замыкается, так как $[\mathcal K]\in H$, а значит является
диаграммой некоторого узла на поверхности $\widetilde S$. Мы можем описать
слабую четность $\psi^H$  и действие соответствующего функториального отображения $\Psi^H$ посредством накрытия $p_H$.

\begin{prop}
\begin{enumerate}
\item перекресток $v$ диаграммы $D$ является четным относительно $\psi^H$ тогда и только тогда, когда
поднятие $\widetilde D_{v}\subset\widetilde D$ половинки $D_{v}$
замыкается;
\item поднятие $\widetilde D$ диаграммы совпадает с результатом применения к диаграмме
$D$ функториального отображения $\Psi^H$.
\end{enumerate}
\end{prop}

\begin{rema}
Во втором пункте утверждения окрестность диаграммы $\widetilde D$  в $\widetilde S$
интерпретируется как поверхность $S'(\widetilde D)$ виртуального
узла (см. рис.~\ref{pic:corvering}).
\begin {figure}[h]
\centering
\includegraphics[height=4cm]{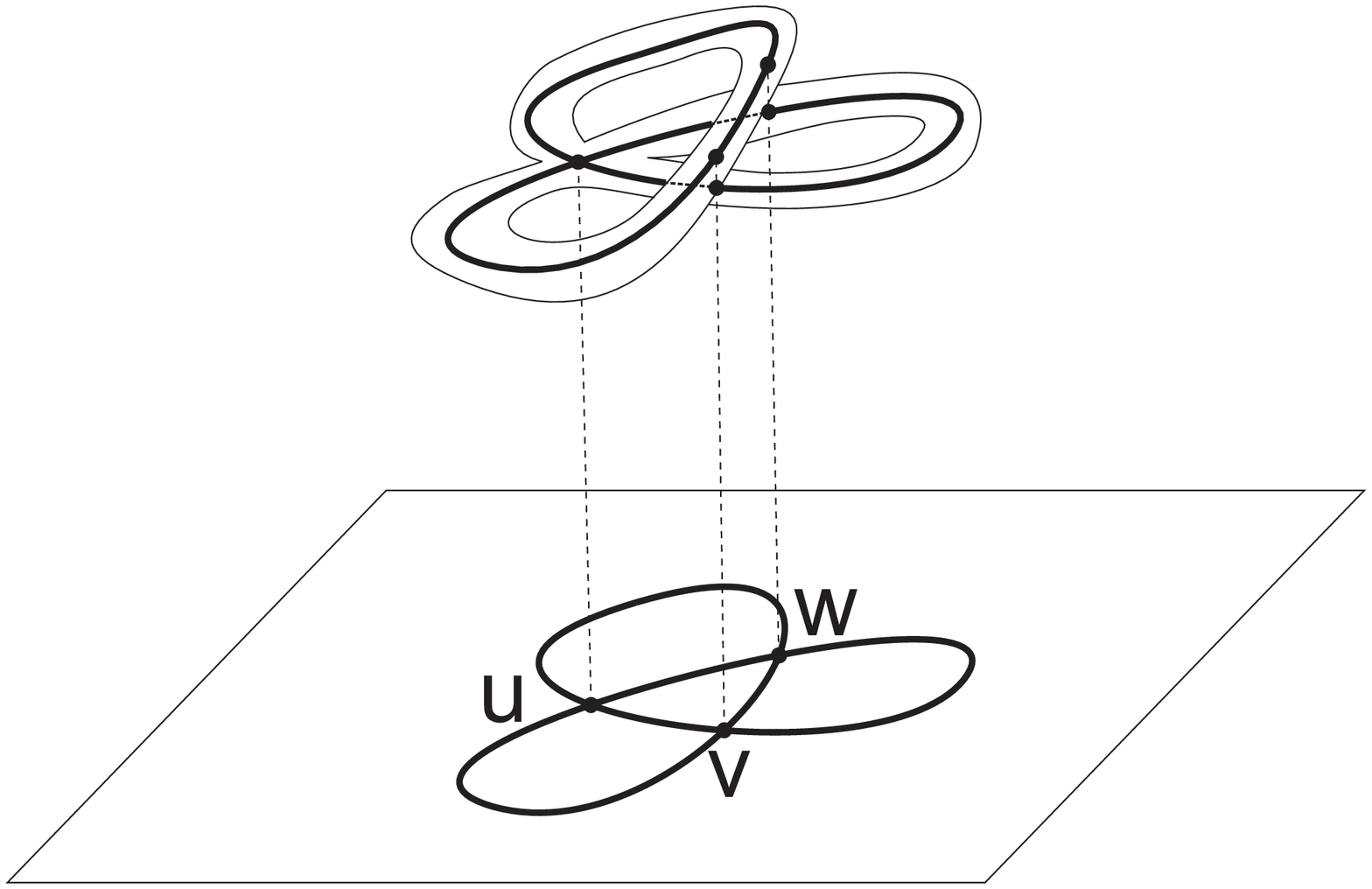}
\caption{Диаграмма и ее поднятие}\label{pic:corvering}
\end {figure}
\end{rema}

\begin{proof}
Для любого перекрестка имеем $\psi^H_D(v)=0$ $\Leftrightarrow$
$[\hat D_v]\in H$ $\Leftrightarrow$ $\widetilde D_{v}$ замыкается.

Второй пункт следует из первого и описания виртуальных диаграмм как
диаграмм в поверхностях.
\end{proof}

В случае гомотопической слабой четности описание функториального
отображения как поднятия при накрытии приводит к следующему
результату.

 \begin{theorem}\label{thm:hom_funct_map}
Пусть $\mathcal K$ --- узел на связной замкнутой ориентированной
поверхности $S$ и $D$  --- его диаграмма. Тогда $\Psi^{hom}(D)$, где
$\Psi^{hom}$ есть функториальное отображение, отвечающее
гомотопической слабой четности, является классической диаграммой
узла.
 \end{theorem}

\begin{proof}
Если $S$ --- сфера, то $D$ --- классическая диаграмма, и отображение
$\Psi^{hom}(D)$ является тождественным.

Пусть род поверхности $S$ больше нуля. Диаграмму $\Psi^{hom}(D)$
можно рассматривать как поднятие диаграммы $D$ на накрывающую
поверхность $\widetilde S$. В случае гомотопической четности
фундаментальная группа $\pi_1(\widetilde S)$ изоморфна  циклической
подгруппе в $\pi_1(\widetilde S)$, порожденной элементом $[\mathcal
K]$. Так как универсальное накрытие поверхности $S$ есть плоскость,
то $\widetilde S$ есть факторизация плоскости по свободно
действующему преобразованию, причем можно считать, что это
преобразование сохраняет метрику постоянной кривизны на плоскости.
Следовательно, поверхность $\widetilde S$ есть плоскость либо
цилиндр. В любом случае, $\widetilde S$, а значит, и окрестность
диаграммы $\Psi^{hom}(D)$ вложима в $\R^2$, так что $\Psi^{hom}(D)$
является классической диаграммой.
\end{proof}

\begin{rema}
Для плоских узлов на фиксированной поверхности
теорема~\ref{thm:hom_funct_map} означает, что функториальное
отображение, соответствующее гомотопической слабой четности,
переводит любой узел в тривиальный, так как таковым является любой
классический плоский узел.
\end{rema}

\begin{corollary}
Мы можем определить проекцию $\Pi$ виртуальных узлов на
классические: для любого виртуального узла $\mathcal K$ положим
$$
\Pi(\mathcal K) = \Psi^{hom}(\overline{\mathcal K}),
$$
где $\overline{\mathcal K}$ есть минимальный представитель узла
$\mathcal K$.
\end{corollary}

\begin{rema}
Построенная проекция определена на уровне узлов, но не диаграмм, то
есть не является функториальным отображением виртуальных узлов в
классические. Отображение (многозначное) на уровне диаграмм из
виртуальных узлов в классические было построено в работе В.О.
Мантурова~\cite{M12}.
\end{rema}


\begin{thebibliography}{99}
%

\bibitem{IMN} D.~Ilyutko, V.~Manturov, I.~Nikonov. Virtual Knot Invariants Arising From Parities. arXiv:1102.5081v1.
\bibitem{IMN_mon}   Ильютко Д.П., Мантуров В.О., Никонов И.М. Четность в теории узлов и граф-зацепления // Современная математика. Фундаментальные направления, т. 41, М.: РУДН, 2011, 163 с.

\bibitem{KK} Kamada N., Kamada S. Abstract link diagrams and virtual knots// J. Knot Theory Ramifications. 2000. 9, \No 1. С. 93--109.
\bibitem{Kuperberg} Kuperberg G. What is a virtual link?// Algebraic and Geometric Topology. 2003. 3. С. 587--591.
\bibitem{M1} 
Мантуров В.О. Четность в теории узлов // Матем. сб., 2010, т. 201, \No 5, С. 65--110. (arXiv:0901.2214v2)
\bibitem{M2} V.~O.~Manturov. On free knots and links/ arXiv:math.GT/0902.0127.
\bibitem{M3} V.~O.~Manturov. Free knots are not invertible/ arXiv:math.GT/0909.2230v2.
\bibitem{M12} V.~Manturov. A Functorial Map from Virtual Knots to Classical Knots and
Generalisations of Parity. arXiv:1011.4640
\bibitem{M4} Мантуров В.О. Четность, свободные узлы, группы и инварианты конечного порядка // Труды ММО, 72:2 (2011),  207-222
\bibitem{M5} 
Мантуров В.О. Четность и кобордизмы свободных узлов // Матем. сб., 203:2 (2012),  45-76. (arXiv:math.GT/1001.2728)
\bibitem{M6} V.~O.~Manturov. Free knots and parity// Introductory lectures on knot theory, 321–345, Ser. Knots Everything, 46, World Sci. Publ., 2012.
\bibitem{M7} Мантуров В.О. Четность и оценка числа виртуальных перекрестков для виртуальных узлов // Труды семинарапо векторному и тензорному анализу, 2012, т. 28, С. 192--210.

\end{thebibliography}
\end{document}